\documentclass[11pt]{amsart}

\usepackage[all]{xy}
\usepackage{color, pstricks}
\usepackage{mathrsfs}
\usepackage{amsmath,amsthm, amssymb, amsgen,amsxtra,amsfonts, amsbsy} 
\usepackage[all]{xy}
\usepackage{verbatim}

\newcommand{\mathds}[1]{{\mathbb #1}}

\begin{document}
%
%   D e f i n i t i o n s
%
%
\theoremstyle{definition}
\newtheorem{Definition}{Definition}[section]
\newtheorem*{Definitionx}{Definition}
\newtheorem{Convention}{Definition}[section]
\newtheorem{Construction}{Construction}[section]
\newtheorem{Example}[Definition]{Example}
\newtheorem{Examples}[Definition]{Examples}
\newtheorem{Exercise}[Definition]{Exercise}
\newtheorem{Remark}[Definition]{Remark}
\newtheorem*{Remarkx}{Remark}
\newtheorem{Remarks}[Definition]{Remarks}
\newtheorem{Caution}[Definition]{Caution}
\newtheorem{Conjecture}[Definition]{Conjecture}
\newtheorem*{Conjecturex}{Conjecture}
\newtheorem{Question}[Definition]{Question}
\newtheorem{Questions}[Definition]{Questions}
\newtheorem*{Acknowledgements}{Acknowledgements}
\newtheorem*{Organization}{Organization}
\newtheorem*{Disclaimer}{Disclaimer}
\theoremstyle{plain}
\newtheorem{Theorem}[Definition]{Theorem}
\newtheorem*{Theoremx}{Theorem}
\newtheorem{Theoremy}{Theorem}
\newtheorem{Proposition}[Definition]{Proposition}
\newtheorem*{Propositionx}{Proposition}
\newtheorem{Lemma}[Definition]{Lemma}
\newtheorem{Corollary}[Definition]{Corollary}
\newtheorem*{Corollaryx}{Corollary}
\newtheorem{Fact}[Definition]{Fact}
\newtheorem{Facts}[Definition]{Facts}
\newtheoremstyle{voiditstyle}{3pt}{3pt}{\itshape}{\parindent}%
{\bfseries}{.}{ }{\thmnote{#3}}%
\theoremstyle{voiditstyle}
\newtheorem*{VoidItalic}{}
\newtheoremstyle{voidromstyle}{3pt}{3pt}{\rm}{\parindent}%
{\bfseries}{.}{ }{\thmnote{#3}}%
\theoremstyle{voidromstyle}
\newtheorem*{VoidRoman}{}

% abgeschrieben aus The LaTeX Companion, 2nd edition,
% von Mittelback & Goossens
%
\newcommand{\prf}{\par\noindent{\sc Proof.}\quad}
\newcommand{\blowup}{\rule[-3mm]{0mm}{0mm}}
\newcommand{\cal}{\mathcal}
\newcommand{\Aff}{{\mathds{A}}}
\newcommand{\BB}{{\mathds{B}}}
\newcommand{\CC}{{\mathds{C}}}
\newcommand{\EE}{{\mathds{E}}}
\newcommand{\FF}{{\mathds{F}}}
\newcommand{\GG}{{\mathds{G}}}
\newcommand{\HH}{{\mathds{H}}}
\newcommand{\NN}{{\mathds{N}}}
\newcommand{\ZZ}{{\mathds{Z}}}
\newcommand{\PP}{{\mathds{P}}}
\newcommand{\QQ}{{\mathds{Q}}}
\newcommand{\RR}{{\mathds{R}}}
\newcommand{\Liea}{{\mathfrak a}}
\newcommand{\Lieb}{{\mathfrak b}}
\newcommand{\Lieg}{{\mathfrak g}}
\newcommand{\Liem}{{\mathfrak m}}
\newcommand{\ideala}{{\mathfrak a}}
\newcommand{\idealb}{{\mathfrak b}}
\newcommand{\idealg}{{\mathfrak g}}
\newcommand{\idealm}{{\mathfrak m}}
\newcommand{\idealp}{{\mathfrak p}}
\newcommand{\idealq}{{\mathfrak q}}
\newcommand{\idealI}{{\cal I}}
\newcommand{\lin}{\sim}
\newcommand{\num}{\equiv}
\newcommand{\dual}{\ast}
\newcommand{\iso}{\cong}
\newcommand{\homeo}{\approx}
\newcommand{\mm}{{\mathfrak m}}
\newcommand{\pp}{{\mathfrak p}}
\newcommand{\qq}{{\mathfrak q}}
\newcommand{\rr}{{\mathfrak r}}
\newcommand{\pP}{{\mathfrak P}}
\newcommand{\qQ}{{\mathfrak Q}}
\newcommand{\rR}{{\mathfrak R}}
%
%  evtl. auch \"uber \mathbb oder \Bbb
%
\newcommand{\OO}{{\cal O}}
\newcommand{\numero}{{n$^{\rm o}\:$}}
\newcommand{\mf}[1]{\mathfrak{#1}}
\newcommand{\mc}[1]{\mathcal{#1}}
\newcommand{\into}{{\hookrightarrow}}
\newcommand{\onto}{{\twoheadrightarrow}}
\newcommand{\Spec}{{\rm Spec}\:}
\newcommand{\BigSpec}{{\rm\bf Spec}\:}
\newcommand{\Spf}{{\rm Spf}\:}
\newcommand{\Proj}{{\rm Proj}\:}
\newcommand{\Pic}{{\rm Pic }}
\newcommand{\Br}{{\rm Br}}
\newcommand{\NS}{{\rm NS}}
\newcommand{\Sym}{{\mathfrak S}}
\newcommand{\Aut}{{\rm Aut}}
\newcommand{\Autp}{{\rm Aut}^p}
\newcommand{\Hom}{{\rm Hom}}
\newcommand{\Ext}{{\rm Ext}}
\newcommand{\ord}{{\rm ord}}
\newcommand{\coker}{{\rm coker}\,}
\newcommand{\divisor}{{\rm div}}
\newcommand{\Def}{{\rm Def}}
\newcommand{\piet}{{\pi_1^{\rm \acute{e}t}}}
\newcommand{\Het}[1]{{H_{\rm \acute{e}t}^{{#1}}}}
\newcommand{\Hfl}[1]{{H_{\rm fl}^{{#1}}}}
\newcommand{\Hcris}[1]{{H_{\rm cris}^{{#1}}}}
\newcommand{\HdR}[1]{{H_{\rm dR}^{{#1}}}}
\newcommand{\hdR}[1]{{h_{\rm dR}^{{#1}}}}
\newcommand{\defin}[1]{{\bf #1}}
\newcommand{\oX}{\cal{X}}
\newcommand{\oA}{\cal{A}}
\newcommand{\oY}{\cal{Y}}
\newcommand{\calC}{{\cal{C}}}
\newcommand{\calL}{{\cal{L}}}
\newcommand{\rmet}{{\rm \acute{e}t}}

\title[Brauer--Severi Varieties and Del~Pezzo Surfaces]{Morphisms to Brauer--Severi varieties, with applications to del~Pezzo surfaces}
\author{Christian Liedtke}
\address{TU M\"unchen, Zentrum Mathematik - M11, Boltzmannstr. 3, D-85748 Garching bei M\"unchen, Germany}
\curraddr{}
\email{liedtke@ma.tum.de}

\date{May 4, 2016}
\subjclass[2010]{14F22,14A10,14J45,14G27}

\begin{abstract}
   We classify morphisms from proper varieties to Brauer--Severi varieties, 
   which generalizes the classical correspondence between morphisms
   to projective space and globally generated invertible sheaves.
   As an application, we study del~Pezzo surfaces of large degree
   with a view towards Brauer--Severi varieties, 
   and recover classical results on rational points, the
   Hasse principle, and weak approximation.
\end{abstract}

\maketitle

\section{Introduction}

\subsection{Overview}
The goal of this article is the study of morphisms $X\to P$ 
from a proper variety $X$ over a field $k$ to a Brauer--Severi variety $P$ over $k$, i.e., $P$
is isomorphic to projective space over the algebraic closure $\overline{k}$ of $k$, 
but not necessarily over $k$.
If $X$ has a $k$-rational point, then so has $P$, and then, $P$ is isomorphic to projective 
space already over $k$.
In this case, there exists a well-known description of morphisms $X\to P$ in terms of
globally generated invertible sheaves on $X$.
However, if $X$ has no $k$-rational point, then we establish in this article
a correspondence between globally generated classes of $\Pic_{(X/k)({\rm fppf})}(k)$,
whose obstruction to coming from an invertible sheaf on $X$ 
is measured by some class $\beta$ in the Brauer group $\Br(k)$, and morphisms to Brauer--Severi
varieties of class $\beta$ over $k$.

As an application of this correspondence, we study del~Pezzo surfaces over $k$ in terms of
Brauer--Severi varieties, and recover many known results
about their geometry and their arithmetic.
If $k$ is a global field, then we obtain applications concerning the Hasse principle and
weak approximation.
Our approach has the advantage of being elementary, self-contained, 
and that we sometimes obtain natural reasons for the existence of $k$-rational points.

\subsection{Morphisms to Brauer--Severi varieties}
Let $X$ be a proper variety over a field $k$, and let $\overline{k}$
be the algebraic closure of $k$.
When studying invertible sheaves on $X$, there are inclusions
and equalities of abelian groups
$$
  \Pic(X) \,\subseteq\, \Pic_{(X/k)(\rmet)}(k)
  \,=\, \Pic_{(X/k)({\rm fppf})}(k)\,\subseteq\, \Pic(X_{\overline{k}}).
$$
On the left (resp. right), we have invertible sheaves on $X$ (resp. $X_{\overline{k}}$)
up to isomorphism, whereas in the middle, we have sections of the 
sheafified relative Picard functor over $k$
(with respect to the \'etale and fppf topology, respectively).
Moreover, the first inclusion is part of an exact sequence
$$
 0\,\to\,\Pic(X) \,\to\, \Pic_{(X/k)(\rmet)}(k)\,\stackrel{\delta}{\longrightarrow}\,\Br(k),
$$
where $\Br(k)$ denotes the Brauer group of the field $k$, and we refer to
Remark \ref{rem: explain delta} for explicit descriptions of $\delta$.
If $X$ has a $k$-rational point, then
$\delta$ is the zero map, i.e., the first inclusion is a bijection.

By definition, a {\em Brauer--Severi variety} is a variety $P$ over $k$,
such that $P_{\overline{k}}\iso\PP_{\overline{k}}^N$ for some $N$, i.e., $P$
is a twisted form of projective space.
Associated to $P$, there exists a Brauer class $[P]\in\Br(k)$ and
by a theorem of Ch\^{a}telet, $P$ is trivial, i.e., 
isomorphic to projective space over $k$, if and only if $[P]=0$.
This is also equivalent to $P$ having a $k$-rational point.
In any case, we have a class $\OO_P(1)\in\Pic_{(P/k)({\rm fppf})}(k)$,
in general not arising from an invertible sheaf on $P$, which becomes
isomorphic to $\OO_{\PP^N}(1)$ over $\overline{k}$,
see Definition \ref{def: O(1) for BS}.

In this article, we extend the notion of a {\em linear system}
to classes in $\Pic_{(X/k)({\rm fppf})}(k)$ that do not necessarily come from
invertible sheaves.
More precisely, we extend the notions of being {\em globally generated}, 
{\em ample}, and {\em very ample} to such classes, see 
Definition \ref{def: globally generated}.
Then, we set up a dictionary between globally generated 
classes in $\Pic_{(X/k)({\rm fppf})}(k)$ and morphisms 
from $X$ to Brauer--Severi varieties over $k$.
In case $X$ has a $k$-rational point, then we  recover the well-known
correspondence between globally generated invertible sheaves and
morphisms to projective space. %, see Remark \ref{rem: trivial case}.
Here is an easy version of our correspondence and we refer to 
Theorem \ref{thm: main} and Remark \ref{rem: trivial case}
for details.

\begin{Theorem}
  \label{theorem1}
  Let $X$ be a proper variety over a field $k$.
   \begin{enumerate}
    \item Let $\varphi:X\to P$ be a morphism to a Brauer--Severi variety $P$ over $k$.
       If we set ${\cal L}:=\varphi^*\OO_P(1)\in\Pic_{(X/k)({\rm fppf})}(k)$,
       then ${\cal L}$ is a globally generated class and
      $$
         \delta({\cal L}) \,=\, [P] \,\in\Br(k).
      $$
    \item If ${\cal L}\in\Pic_{(X/k)({\rm fppf})}(k)$ is globally generated,
     then $\calL\otimes_k\overline{k}$ corresponds to a unique invertible
     sheaf ${\cal M}$ on $X_{\overline{k}}$ and the morphism associated 
     to the complete linear system $|{\cal M}|$ 
     descends to a morphism over $k$
     $$
        |{\cal L}|\,:\,X \,\to\, P,    
     $$
     where $P$ is a Brauer--Severi variety over $k$ with
     $\delta({\cal L})=[P]$.
  \end{enumerate}
\end{Theorem}

We note that our result is inspired by a geometric construction of Brauer--Severi varieties
of Grothendieck, see \cite[Section (5.4)]{Grothendieck Brauer}, and it seems that
it is known to the experts.
As immediate corollaries, we recover two classical theorems about Brauer--Severi
varieties due to Ch\^{a}telet and Kang,
see Corollary \ref{cor: kang} and Corollary \ref{cor: chatelet}.
\medskip

\subsection{Del~Pezzo surfaces}
In the second part, we apply this machinery to the geometry and arithmetic of 
del~Pezzo surfaces over arbitrary ground fields.
I would like to stress that most, if not all, of the results of this second part
are well-known. 
To the best of my knowledge, I have tried to give the original references.
However, my organization of the material and the hopefully more geometric
approach to del~Pezzo surfaces via morphisms to Brauer--Severi varieties is new.

By definition, a {\em del~Pezzo surface} is a smooth and proper surface
$X$ over a field $k$, whose anti-canonical invertible sheaf $\omega_X^{-1}$ is ample.
The {\em degree} of a del~Pezzo surface is the self-intersection 
number of $\omega_X$.
The classification of del~Pezzo surfaces over $\overline{k}$ is well-known:
The degree $d$ satisfies $1\le d\le 9$, and they are isomorphic either 
to $\PP^1\times\PP^1$ or 
to the blow-up of $\PP^2$ in $(9-d)$ points in general position.

As an application of Theorem \ref{theorem1}, % (= Theorem \ref{thm: main}), 
we obtain the following.
\begin{enumerate}
 \item If $d=8$ and $X_{\overline{k}}\iso\PP^1_{\overline{k}}\times\PP^1_{\overline{k}}$,
   then there exists an embedding
   $$\begin{array}{ccccc}
      |-\frac{1}{2}K_X| &:& X &\into& P
       \end{array}
   $$
   into a Brauer--Severi threefold $P$.
   Moreover, $X$ is either isomorphic to a product  of two
   Brauer--Severi curves or to a quadratic twist
   of the self-product of a Brauer--Severi curve.
   We refer to Theorem \ref{thm: product type} and 
   Proposition \ref{prop: product type classification} for details.
 \item If $d\geq7$ and $X_{\overline{k}}\not\iso\PP^1_{\overline{k}}\times\PP^1_{\overline{k}}$,
   then there exists a birational morphism
   $$\begin{array}{ccccc}
      f&:&X&\to&P
      \end{array}
   $$
   to a Brauer--Severi surface $P$ over $k$ that is the blow-up in a closed and zero-dimensional
   subscheme of length $(9-d)$ over $k$.
   We refer to Theorem \ref{thm: del Pezzo descent} for details.
  \item If $d=6$, then there exist two finite field extensions $k\subseteq K$ and $k\subseteq L$
   with $[K:k]|2$ and $[L:k]|3$ such that there exists a birational morphism $f:X\to P$
   to a Brauer--Severi surface $P$ over $k$ that is the blow-up in a closed and zero-dimensional
   subscheme of length $3$ over $k$ if and only $k=K$.
   On the other hand, there exists a birational morphism $X\to Y$ onto a degree 
   $8$ del~Pezzo surface $Y$ 
   of product type if and only if $k=L$.
   We refer to Theorem \ref{thm: degree 6} for details.
 \item For partial results if $d\leq5$, as well as birationality criteria for
   when a del~Pezzo surface
   is birationally equivalent to a Brauer--Severi surface, we refer to Section \ref{sec: small degree}. 
\end{enumerate}

As further applications, we recover well-known results about rationality, 
unirationality, existence of $k$-rational points, Galois cohomology, 
the Hasse principle, and weak approximation for del~Pezzo surfaces.

\begin{Acknowledgements}
It is a pleasure for me to thank J\"org Jahnel, Andrew Kresch, Raphael Riedl, 
Ronald van~Luijk, and Anthony V\'arilly-Alvarado for comments and discussions.
I especially thank Jean-Louis Colliot-Th\'el\`ene and Alexei Skorobogatov
for providing me with references,
discussions, and pointing out mistakes, as well as correcting some of my too naive ideas.
Last, but not least, I thank the referee for careful proof-reading and the many useful
suggestions.
\end{Acknowledgements}

\section*{Notations and Conventions}

In this article, $k$ denotes an arbitrary field,
$\overline{k}$ (resp. $k^{\rm sep}$) its algebraic (resp. separable) closure, and 
$G_k={\rm Gal}(k^{\rm sep}/k)$ its absolute Galois group.
By a variety over $k$ we mean a scheme $X$ that is of finite type, separated,
and geometrically integral over $k$.
If $K$ is a field extension of $k$, then we define $X_K:=X\times_{\Spec k}\Spec K$.

\section{Picard functors and Brauer groups}
\label{sec:brauer}

This section, we recall a couple of definitions and general results
about the various relative Picard functors, about Brauer groups of fields and schemes,
as well as Brauer--Severi varieties.

\subsection{Relative Picard functors}
Let us first recall a couple of generalities about the several Picard functors.
Our main references are \cite{Grothendieck Picard1}, \cite{Grothendieck Picard2}, 
as well as the surveys \cite[Chapter 8]{BLR} and \cite{Kleiman Picard}.

For a scheme $X$, we define its {\em Picard group} $\Pic(X)$ to be the abelian group of 
invertible sheaves on $X$ modulo isomorphism.
If $f:X\to S$ is a separated morphism of finite type over a Noetherian base scheme $S$,
then we define the {\em absolute Picard functor} to be the functor that associates to 
each Noetherian $T\to S$ the abelian group $\Pic_X(T):=\Pic(X_T)$, where $X_T:=X\times_ST$.
Now, as explained, for example in \cite[Section 9.2]{Kleiman Picard}, the absolute Picard functor is 
a separated presheaf for the Zariski, \'etale, and the fppf topologies, but it is never a sheaf for 
the Zariski topology.
In particular, the absolute Picard functor is never representable by a scheme or by an algebraic space.
This leads to the introduction of the {\em relative Picard functor} $\Pic_{X/S}$ by setting
$\Pic_{X/S}(T):=\Pic(X_T)/\Pic(T)$, and then, we have the associated sheaves for the Zariski, \'etale, 
and fppf topologies
$$
  \Pic_{(X/S)({\rm zar})},\mbox{ \quad }   \Pic_{(X/S)({\rm \acute{e}t})},\mbox{ \quad and \quad }   \Pic_{(X/S)({\rm fppf})}.
$$
In many important cases, these sheaves are representable by schemes or algebraic spaces over $S$.
For our purposes, it suffices to work with the sheaves so that we will not address representability questions here, 
but refer the interested reader to \cite[Chapter 8.2]{BLR} and \cite[Chapter 9.4]{Kleiman Picard} instead.
Having introduced these sheaves, let us recall the following easy facts, 
see, for example, \cite[Exercise 9.2.3]{Kleiman Picard}.

\begin{Proposition}
  \label{easy picard facts}
  Let $X\to S$ be a scheme that is separated and of finite type over a Noetherian scheme $S$.
  Let $L$ be a field with a morphism $\Spec L\to S$.
  \begin{enumerate}
   \item Then, the following natural maps are isomorphisms:
  $$
     \Pic_X(L)\,\stackrel{\iso}{\longrightarrow}\,\Pic_{X/S}(L)\,\stackrel{\iso}{\longrightarrow}\,\Pic_{(X/S)({\rm zar})}(L).
  $$
   \item If $L$ is algebraically closed, then also the following natural maps are isomorphisms:
  $$
    \Pic_X(L)\,\stackrel{\iso}{\longrightarrow}\,\Pic_{(X/S)({\rm \acute{e}t})}(L)\,\stackrel{\iso}{\longrightarrow}\,\Pic_{(X/S)({\rm fppf})}(L).
  $$
  \end{enumerate}
 \end{Proposition}
 
It is important to note that if $L$ is not algebraically closed, then the
natural map $\Pic_X(L)\to\Pic_{(X/S)({\rm \acute{e}t})}(L)$ is usually not an isomorphism, 
i.e., not every section of $\Pic_{(X/S)({\rm \acute{e}t})}$ over $L$
arises from an invertible sheaf on $X_L$.
The following example, taken from \cite[Exercise 9.2.4]{Kleiman Picard}, 
is crucial to everything that follows and illustrates this.

\begin{Example}
  \label{ex: failure picard}
  Let $X$ be the smooth plane conic over $\RR$ defined by
  $$
     X\,:=\,\{\, x_0^2+x_1^2+x_2^2=0 \,\} \,\subset\,\PP^2_{\RR}.
  $$
  Then, $X$ is not isomorphic to $\PP^1_\RR$ since $X(\RR)=\emptyset$, 
  but there exists an isomorphism $X_\CC\to\PP^1_\CC$.
  In particular, $X$ is an example of a non-trivial Brauer--Severi variety (see Definition \ref{def: brauer severi}).
  
  Next, if $x\in X$ is a closed point, then $\kappa(x)\iso\CC$, that is, $x$ is a zero-cycle of degree $2$.
  Moreover, $\OO_X(x)$ generates $\Pic_X(\RR)$, for if there was an invertible sheaf of
  odd degree on $X$, then there would exist an invertible sheaf of degree $1$ on $X$ and then,
  Riemann--Roch would imply $X(\RR)\neq\emptyset$, a contradiction. 
  
  On the other hand, $x$ splits on $X_\CC$ into two closed points, say $x_1$ and $x_2$.
  Since $\OO_{X_\CC}(x_1)$ and $\OO_{X_\CC}(x_2)$ are isomorphic as invertible sheaves on $X_\CC$, 
  it follows that $\OO_{X_\CC}(x_1)$ descends from a class in $\Pic_{(X/\RR)({\rm \acute{e}t})}(\CC)$
  to a class in $\Pic_{(X/\RR)({\rm \acute{e}t})}(\RR)$.
  
  These observations show that the natural map
  $\Pic_X(\RR)\to\Pic_{(X/\RR)({\rm \acute{e}t})}(\RR)$ is not surjective.
\end{Example}

In this example, we have $X(\RR)=\emptyset$, i.e., the structure morphism $X\to\Spec\RR$ has no section.
Quite generally, we have the following comparison theorem for the several relative Picard functors,
and refer, for example, to \cite[Theorem 9.2.5]{Kleiman Picard} for details and proofs.

\begin{Theorem}[Grothendieck] \label{thm: comparison}
  Let $f:X\to S$ be a scheme that is separated and of finite type over a Noetherian scheme $S$,
  and assume that $\OO_S\stackrel{\iso}{\longrightarrow}f_\ast\OO_X$ holds universally.
  \begin{enumerate}
    \item Then, the natural maps 
    $$
       \Pic_{X/S} \,\into\, \Pic_{(X/S)({\rm zar})} \,\into\,  \Pic_{(X/S)({\rm \acute{e}t})} \,\into\,   \Pic_{(X/S)({\rm fppf})}
    $$
    are injections.
    \item If $f$ has a section, then all three maps are isomorphisms.
     If $f$ has a section locally in the Zariski topology, then the latter two maps are isomorphisms, and if
     $f$ has a section locally in the \'etale topology, then the last map is an isomorphism.
  \end{enumerate}
\end{Theorem}

To understand the obstruction to realizing a section of
$\Pic_{(X/S)(\rmet)}$ or $\Pic_{(X/S)({\rm fppf})}$ over $S$
by an invertible sheaf on $X$ in case there is no section of $X\to S$, 
we recall the following definition.

\begin{Definition}
  \label{def: brauer}
  For a scheme $T$, the \'etale cohomology group $\Het{2}(T,\GG_m)$ is called the 
  {\em cohomological Brauer group}, and is denoted $\Br'(T)$.
  The set of sheaves of Azumaya algebras on $T$ modulo Brauer equivalence 
  also forms a group, the {\em Brauer group} of $T$, and is denoted $\Br(T)$.
\end{Definition}

We will not discuss sheaves of Azumaya algebras on schemes in the sequel, 
but only remark that these generalize central simple algebras over fields 
(see Section \ref{subsec: BS} for the latter), and refer the interested reader to
\cite{Grothendieck Brauer1} and \cite[Chapter IV]{Milne} for details and references, 
as well as to \cite{Poonen} for a survey.

Using that $\GG_m$ is a smooth group scheme, Grothendieck \cite{Grothendieck Brauer} showed that
the natural map $\Het{2}(T,\GG_m)\to H^2_{{\rm fppf}}(T,\GG_m)$ is an isomorphism, i.e., 
it does not matter whether the cohomological Brauer group $\Br'(T)$ is defined with respect to
the \'etale or the fppf topology.
Next, there exists a natural injective group homomorphism $\Br(T)\to\Br'(T)$, 
whose image is contained in the torsion subgroup of $\Br'(T)$.
If $T$ is the spectrum of a field $k$, then this injection is even an isomorphism, i.e., 
$\Br(k)=\Br'(k)$, see, for example, \cite{Grothendieck Brauer}, \cite{Gille Szamuely}, and 
\cite[Chapter IV]{Milne} for details and references.

The connection between Brauer groups, Proposition \ref{easy picard facts}, 
and Theorem \ref{thm: comparison} is as follows, see, 
for example \cite[Chapter 8.1]{BLR} or \cite[Section 9.2]{Kleiman Picard}.

\begin{Proposition}
  \label{prop: delta}
   Let $f:X\to S$ be a scheme that is separated and of finite type over a Noetherian scheme $S$,
  and assume that $\OO_S\stackrel{\iso}{\longrightarrow}f_\ast\OO_X$ holds universally.
  Then, for each $S$-scheme $T$
  there exists a canonical exact sequence
  $$
    0\,\to\,\Pic(T)\,\to\,\Pic(X_T)\,\to\,\Pic_{(X/S)({\rm fppf})}(T)\,\stackrel{\delta}{\longrightarrow}\,\Br'(T)\,\to\,\Br'(X_T)\,.
  $$
  If $f$ has a section, then $\delta$ is the zero-map.\qed
\end{Proposition}

\subsection{Varieties and the Amitsur subgroup}
By our conventions above, a variety over a field $k$ is
a scheme $X$ that is of finite type, separated,
and geometrically integral over $k$.
In this situation, the conditions of Proposition \ref{prop: delta} are fulfilled,
as the following remark shows.

\begin{Remark}
 \label{rem: geometry}
 If $X$ is a proper variety over a field $k$, then
 \begin{enumerate}
  \item the structure morphism $f:X\to \Spec k$ is separated, of finite type, and 
 $\OO_{\Spec k}\iso f_\ast\OO_X$ holds universally.
 \item The morphism $f$ has sections locally in the \'etale topology (see, for example,
   \cite[Appendix A]{Gille Szamuely}).
 \item Since the base scheme is a field $k$, we have $\Br(k)=\Br'(k)$.
 \end{enumerate}
 In Remark \ref{rem: explain delta}, we will give an explicit description of $\delta$ in this case.
\end{Remark}

In Example \ref{ex: failure picard}, the obstruction to representing the
class of ${\cal L}:=\varphi^*\OO_{\PP^1_\CC}(1)$ in $\Pic_{(X/\RR)({\rm fppf})}(\RR)$ by an invertible sheaf
on $X$ can be explained via $\delta$, which maps $\cal L$ to the non-zero element of
$\Br(\RR)\iso\ZZ/2\ZZ$.
In terms of Azumaya algebras (since the base is $\Spec\RR$, these are
central simple $\RR$-algebras),
this Brauer class corresponds the $\RR$-algebra $\HH$ of quaternions, but we will not 
pursue this point of view in the sequel.

\begin{Proposition}
  \label{prop: picard in geometry}
  Let $X$ be a proper variety over a field $k$.
  Then, there exist natural isomorphisms of abelian groups
  $$
    \Pic_{X/k}(k^{\rm sep})^{G_k} \,\stackrel{\iso}{\longrightarrow}\,
    \Pic_{(X/k)(\rmet)}(k) \,\stackrel{\iso}{\longrightarrow}\,   \Pic_{(X/k)({\rm fppf})}(k),
  $$
  where the ${-}^{G_k}$ denotes Galois invariants.
\end{Proposition}

\prf
The first isomorphism follows from Galois theory and sheaf axioms and the second isomorphism
follows from Theorem \ref{thm: comparison} and Remark \ref{rem: geometry}.
\qed\medskip

The Brauer group $\Br(k)$ of a field $k$ is an abelian torsion group, see,
for example, \cite[Corollary 4.4.8]{Gille Szamuely}.
Motivated by Proposition \ref{prop: delta}, we introduce the following
subgroup of $\Br(k)$ that measures the deviation between
$\Pic_{(X/k)({\rm fppf})}(k)$ and $\Pic(X)$.

\begin{Definition}
 \label{def: amitsur group}  
   Let $X$ be a proper variety over a field $k$.
  Then, the {\em Amitsur subgroup} of $X$ in $\Br(k)$ 
  is the subgroup 
  $$
     {\rm Am}(X) \,:=\, \delta(\Pic_{(X/k)({\rm fppf})}(k))\,\subseteq\,\Br(k).
  $$
  By the previous remarks, it is an abelian torsion group.
\end{Definition}

The following lemma gives bounds for the order of torsion in ${\rm Am}(X)$.

\begin{Lemma}
  \label{lem: order amitsur}
  Let $X$ be a proper variety over a field $k$.
  If there exists a closed point on $X$, whose residue field is
  of degree $n$ over $k$, then every element of ${\rm Am}(X)$ has an order
  dividing $n$.
\end{Lemma}

\prf
Let  $x\in X$ be a closed point, say, with residue field $K/k$ that is of degree $n$
over $k$.
Since $X_K$ has a $K$-rational point, the map $\delta$ of $X_K$ is identically
zero by Proposition \ref{prop: delta}.
Thus, we have an inclusion ${\rm Am}(X) \subseteq\Br(K|k):=\ker(\Br(k)\to\Br(K))$,
where $\Br(k)\to\Br(K)$ is the restriction homomorphism.

If $K$ is separable over $k$, then $\Br(K|k)$ is contained in the
$n$-torsion of $\Br(k)$, which follows from the fact that
the composition of restriction and corestriction
is multiplication by $n$, see \cite[Proposition 4.2.10]{Gille Szamuely}.

If $K$ is a purely inseparable extension of $k$, generated by $p^r$-th roots,
then $\Br(K|k)$ is $p^r$-torsion (which yields even stronger bounds on the
torsion than claimed), see for example, Hochschild's Theorem
\cite[Theorem 9.1.1]{Gille Szamuely} for an explicit description for
this group.

In general, we can factor the extension $K/k$ into a separable and a purely inseparable 
extension, and by combining the previous two special cases, the statement follows.
\qed\medskip

Using Proposition \ref{prop: delta}, we can give two
alternative definitions of ${\rm Am}(X)$.
In fact, the birational invariance of this group for
Brauer--Severi varieties is a classical result of Amitsur, probably known 
to Ch\^{a}telet  and Witt in some form or another,
see also Theorem \ref{thm: amitsur birational} below.

\begin{Proposition}
 \label{prop: amitsur is birational invariant}
 Let $X$ be a smooth and proper variety over $k$.
 Then,
 $$
   {\rm Am}(X)\,=\, \ker\left(\Br(k)\,\to\,\Br'(X)\right)
   \,=\,\ker\left(\Br(k)\,\to\,\Br(k(X))\right).
 $$
 In particular, ${\rm Am}(X)$ is a birational invariant
 of smooth and proper varieties over $k$.
\end{Proposition}

\prf
The first equality follows from the exact sequence of
Proposition \ref{prop: delta}.
Since $X$ is smooth over $k$, 
the natural map $\Br'(X)\to\Br(k(X))$ is injective,
see, for example, \cite[Example III.2.22]{Milne},
and then, the second equality follows.
From this last description, it is clear that ${\rm Am}(X)$
is a birational invariant.
\qed\medskip

\begin{Remark}
  In \cite[Section 5]{CTKMdP6}, the kernel of $\Br(k)\to\Br(k(X))$ was denoted 
  $\Br(k(X)/k)$.
  Thus, if $X$ is smooth and proper over $k$, then this latter group
  coincides with ${\rm Am}(X)$.
 However, this group should not be confused with $\Br(k(X))/\Br(k)$,
 which is related to another important birational invariant that we will introduce in
 Section \ref{subsec: dP arithmetic}.
\end{Remark}

If $X$ has a $k$-rational point, then ${\rm Am}(X)=0$ by
Proposition \ref{prop: delta}.
On the other hand, there exist proper varieties $X$ with trivial
Amitsur subgroup without $k$-rational points (some
degree $8$ del~Pezzo surfaces of product type with $\rho=1$
provide examples, see Proposition \ref{cor: h1 product type}).
Let us recall that a {\em zero-cycle} on $X$ is a formal finite sum 
$\sum_i n_i Z_i$, where the $n_i\in\ZZ$ and where the $Z_i$ are closed
points of $X$.
It is called {\em effective} if $n_i\geq0$ for all $i$.
The {\em degree} is defined to be $\deg(Z):=\sum_i n_i[\kappa(Z_i):k]$,
where $\kappa(Z_i)$ denotes the residue field of the point $Z_i$.

\begin{Corollary}
 \label{cor: amitsur trivial}
  Let $X$ be a proper variety over a field $k$.
  If there exists a zero cycle of degree $1$ 
  on $X$, then ${\rm Am}(X)=0$.
 \qed
\end{Corollary} 

If $X$ is a projective variety over $k$, then
$\Pic_{(X/k)(\rmet)}$ and $\Pic_{(X/k)({\rm fppf})}$ are representable by a group scheme
$\Pic_{X/k}$ over $k$, the {\em Picard scheme}.
The connected component of the identity is denoted $\Pic_{X/k}^0$, and
the quotient 
$$
  \NS_{X/k}(\overline{k}) \,:=\, \Pic_{X_{\overline{k}}/\overline{k}}(\overline{k}) \,/\, \Pic^0_{X_{\overline{k}}/\overline{k}}(\overline{k}),
$$
the {\em N\'eron--Severi group}, is a finitely generated abelian group, whose rank is
denoted $\rho(X_{\overline{k}})$.
We refer to \cite[Section 8.4]{BLR} for further discussion.
Moreover, if $X$ is smooth over $k$, then $\Pic_{X/k}^0$ is of dimension 
$\frac{1}{2}b_1(X)$, where $b_1$ denotes the first $\ell$-adic Betti number.

\begin{Lemma}
  \label{lem: picard rank}
  Let $X$ be a smooth and projective variety over a field $k$ with $b_1(X)=0$.
  Then, $\Pic_{(X/k)({\rm fppf})}(k)$ is a finitely generated abelian group,
  $$
    {\rm rank}\,\Pic(X)\,=\,{\rm rank}\,\Pic_{(X/k)({\rm fppf})}(k) \,\leq\,\rho(X_{\overline{k}}),
  $$
  and ${\rm Am}(X)$ is a finite abelian group.
\end{Lemma}

\prf
If $b_1(X)=0$, then, by the previous discussion, $\Pic(X_{\overline{k}})$ is a finitely generated
abelian group of rank $\rho(X_{\overline{k}})$.
Since $\Pic(X)$ and $\Pic_{(X/k)({\rm fppf})}(k)$ are contained in $\Pic(X_{\overline{k}})$, they
are also finitely generated of rank at most $\rho(X_{\overline{k}})$.
Since ${\rm Am}(X)=\delta(\Pic_{(X/k)({\rm fppf})}(k))$ is a torsion subgroup of $\Br(k)$,
Proposition \ref{prop: delta} implies the stated equality of ranks.
Moreover, being torsion and a finitely generated abelian group, ${\rm Am}(X)$ is finite.
\qed\medskip

\subsection{Brauer--Severi varieties}
\label{subsec: BS}
Next, we recall a couple of results about Brauer--Severi varieties, and refer the interested reader
to \cite[Chapter 5]{Gille Szamuely} and the surveys \cite{Jahnel}, \cite{Poonen} 
for details, proofs, and further references.

\begin{Definition}
  \label{def: brauer severi}
  A {\em Brauer--Severi variety} over a field $k$ is a proper variety $P$ over $k$, 
  such that there exists a finite field extension $K$ of $k$ and an
  isomorphism $P_K\iso \PP_K^n$ over $K$.
\end{Definition} 

In case $P$ is of dimension one (resp. two, resp. three), we will also
refer to it as a Brauer--Severi curve (resp. Brauer--Severi surface, resp. Brauer--Severi threefold).
Any field extension $K$ of $k$ such that $P_K$ is isomorphic to projective space over $K$ is called
a {\em splitting field} for $P$, and $P$ is said to {\em split} over $K$.
By a theorem of Ch\^{a}telet, a Brauer--Severi variety $P$ over $k$ is {\em trivial}, i.e.,  splits over $k$,
i.e., is $k$-isomorphic to projective space over $k$, if and only if it possesses a $k$-rational point.
Since a geometrically integral variety over a field $k$ always has points over $k^{\rm sep}$,
it follows that a Brauer--Severi variety can be split over a finite and separable extension of $k$, which 
we may also assume to be Galois if we want.

For a finite field extension $K$ of $k$ that is Galois with Galois group $G$, the set of all Brauer--Severi varieties
of dimension $n$ over $k$ that split over $K$, can be interpreted as the set of all $G$-twisted forms of $\PP^n_K$, 
which is in bijection to the cohomology group $H^1(G,{\rm Aut}(\PP^n_K))$.
Using ${\rm Aut}(\PP^n)\iso{\rm PGL}_{n+1}$, and taking cohomology in the short
exact sequence
$$
    1\,\to\,\GG_m\,\to\,{\rm GL}_{n+1}\,\to\,{\rm PGL}_{n+1}\,\to\,1,
$$
the boundary map associates to the class of a Brauer--Severi variety $P$ of dimension
$n$ in $H^1(G,{\rm PGL}_{n+1}(K))$ a class in 
$$
  \Br(K|k) \,:=\, \ker\left(\Br(k)\to\Br(K)\right) \,=\,\ker\left( \Het{2}(k,\GG_m) \to \Het{2}(K,\GG_m) \right).
$$
Taking the limit over all finite Galois extensions of $k$, we obtain for every Brauer--Severi
variety $P$ over $k$ a class $[P]\in\Br(k)$.
This cohomology class is torsion and its order is called the {\em period} of $P$, denoted
${\rm per}(P)$.
By a theorem of Ch\^{a}telet, a Brauer--Severi variety is trivial if and only if the class
$[P]\in\Br(k)$ is zero, i.e., if and only if ${\rm per}(P)=1$.
We will say that two Brauer--Severi varieties over $k$ are {\em Brauer equivalent}
if their associated classes in $\Br(k)$ are the same.

To say more about Brauer classes associated to Brauer--Severi varieties, 
we will shortly digress on non-commutative $k$-algebras,
and refer to \cite[Section 2]{Gille Szamuely} and \cite{Jacobson} for details:
We recall that a {\em central simple $k$-algebra} is a 
$k$-algebra $A$, whose center is equal to $k$ (i.e., $A$ is central), and 
whose only two-sided ideals are $(0)$ and $A$ (i.e., $A$ is simple).
If $A$ is moreover finite-dimensional over $k$, then by theorems of Noether, K\"othe, and Wedderburn, 
there exists a finite and separable field extension $k\subseteq K$ that {\em splits} $A$, i.e.,
$A\otimes_kK\iso{\rm Mat}_{n\times n}(K)$.
In particular, the dimension of $A$ over $k$ is always a square,
and we set the {\em degree} of $A$ to be
${\rm deg}(A):=\sqrt{\dim_k(A)}$.
Two central simple $k$-algebras $A_1$ and $A_2$ are said to be {\em Brauer equivalent}
if there exist integers $a_1,a_2\geq1$ such that
$A_1\otimes_k{\rm Mat}_{a_1\times a_1}(k)\iso A_2\otimes_k{\rm Mat}_{a_2\times a_2}(k)$.

The connection between central simple algebras and Brauer--Severi varieties
is the following dictionary, see \cite[Theorem 2.4.3]{Gille Szamuely}.

\begin{Theorem}
  \label{thm: central simple algebras}
  Let $k\subseteq K$ be a field extension that is Galois with Galois group $G$.
  Then, there is a natural bijection of sets between
  \begin{enumerate}
    \item Brauer--Severi varieties of dimension $n$ over $k$ that split over $K$,
    \item $H^1(G, {\rm PGL}_{n+1}(K))$, and
    \item central simple $k$-algebras of  degree $n+1$ over $k$ that split over $K$.
  \end{enumerate}
  Under this bijection, Brauer equivalence of (1) and (3) coincide.
\end{Theorem}

We also recall that a {\em division algebra} is a $k$-algebra in which
every non-zero element has a two-sided multiplicative inverse.
For example, field extensions of $k$ are division algebras, and a non-commutative
example is provided by the quaternions over $\RR$.
Given a simple and finite-dimensional $k$-algebra $A$,
a theorem of Wedderburn states that there exists a unique division algebra $D$
over $k$ and a unique integer $m\geq1$
and an isomorphism of $k$-algebras
$A\iso{\rm Mat}_{m\times m}(D)$, see \cite[Theorem 2.1.3]{Gille Szamuely}.

\begin{Corollary}
  \label{cor: isomorphic BS}
  If two Brauer--Severi varieties over $k$ of the
  same dimension are Brauer equivalent,
   then they are isomorphic as schemes over $k$.
\end{Corollary}

\prf
By Theorem \ref{thm: central simple algebras}, it suffices to show that
two Brauer equivalent central simple $k$-algebras $A_1$, $A_2$ of the same dimension
are isomorphic.
By Wedderburn's theorem, there exist
division algebras $D_i$ and integers $m_i\geq1$ such that
$A_i\iso{\rm Mat}_{m_i\times m_i}(D_i)$ for $i=1,2$.
By definition of Brauer-equivalence, there exist integers $a_i\geq1$ and
an isomorphism of $k$-algebras
$$
 A_1\otimes_k{\rm Mat}_{a_1\times a_1}(k) \,\iso\, A_2\otimes_k{\rm Mat}_{a_2\times a_2}(k).
$$
Together with the $k$-algebras isomorphisms
$$
 \begin{array}{lcl}
  A_i\otimes_k{\rm Mat}_{a_i\times a_i}(k)  
  &\iso&  {\rm Mat}_{m_i\times m_i}(D_i)\otimes_k{\rm Mat}_{a_1\times a_1}(k) \\
  &\iso&  {\rm Mat}_{a_im_i\times a_im_i}(D_i)
  \end{array}
$$
and the uniqueness part in Wedderburn's theorem, we conclude
$D_1\iso D_2$, as well as $a_1=a_2$, whence $A_1\iso A_2$,
see also \cite[Remark 2.4.7]{Gille Szamuely}.
\qed\medskip

For Brauer--Severi varieties over $k$ that are of different dimension,
we refer to Ch\^{a}telet's theorem (Corollary \ref{cor: chatelet}) below.
On the other hand, for Brauer--Severi varieties over $k$
that are of the same dimension, 
Amitsur conjectured that they are birationally equivalent if and only if their classes
generate the same cyclic subgroup of $\Br(k)$, see also
Remark \ref{rem: the amitsur conjecture}.

For projective space, the degree map  $\deg:\Pic(\PP^n_k)\to\ZZ$, which sends
$\OO_{\PP^n_k}(1)$ to $1$, is an isomorphism.
Thus, if $P$ is a Brauer--Severi variety over $k$ and $G_k:={\rm Gal}(k^{\rm sep}/k)$,
then there are isomorphisms
$$\begin{array}{lclcl}
 \Pic_{(P/k)({\rm fppf})}(k) &\iso&  \Pic_{(P/k)}(k^{\rm sep})^{G_k} &\iso&\Pic_{(P/k)}(k^{\rm sep}) \\
   &\iso&  \Pic(\PP^{\dim(P)}_{k^{\rm sep}}) &\stackrel{{\rm deg}}{\longrightarrow}& \ZZ.
\end{array}$$
The first isomorphism is Proposition \ref{prop: picard in geometry}, and the second follows
from the fact that the $G_k$-action must send the unique ample generator of $\Pic_{(P/k)}(k^{\rm sep})$
to an ample generator, showing that $G_k$ acts trivially.
The third isomorphism follows from the fact that $P$ splits over 
a separable extension.

\begin{Definition}
   \label{def: O(1) for BS}
   For a Brauer--Severi variety $P$ over $k$, we denote the
   unique ample generator of $\Pic_{(P/k)({\rm fppf})}(k)$
   by $\OO_P(1)$.
\end{Definition}

We stress that $\OO_P(1)$ is a class in $\Pic_{(P/k)({\rm fppf})}(k)$ that 
usually does not come from an invertible sheaf on $P$ - in fact this happens
if and only if $P$ is a trivial Brauer--Severi variety, i.e., split over $k$.
For a Brauer--Severi variety, the short exact sequence from Proposition \ref{prop: delta}
becomes the following.

\begin{Theorem}[Lichtenbaum]
  \label{thm: brauer severi picard}
  Let $P$ be a Brauer--Severi variety over $k$. Then, there exists an exact sequence
  $$
   0\,\to\, \Pic(P) \,\to\, 
     \underbrace{\Pic_{(P/k)({\rm fppf})}(k)}_{\iso\,\ZZ} 
   \,\stackrel{\delta}{\longrightarrow}\, \Br(k) 
   \,\to\,\Br(k(P))\,.
  $$
  More precisely, we have
  \begin{eqnarray*}
     \delta(\OO_P(1)) &=& [P], \mbox{ \quad and} \\
     \Pic(P)        &=& \OO_P({\rm per}(P))\cdot\ZZ.
   \end{eqnarray*}
  Since $\omega_P\iso\OO_P(-\dim(P)-1)$, the period ${\rm per}(P)$
  divides $\dim(P)+1$.
\end{Theorem}

Again, we refer to \cite[Theorem 5.4.5]{Gille Szamuely} for details and proofs.
Using Proposition \ref{prop: amitsur is birational invariant}, we immediately
obtain the following classical result of Amitsur \cite{Amitsur} as corollary.

\begin{Theorem}[Amitsur]
  \label{thm: amitsur birational}
  If $P$ is a Brauer--Severi variety over $k$, then ${\rm Am}(P)\iso\ZZ/{\rm per}(P)\ZZ$.
  If two Brauer--Severi varieties are birationally equivalent over $k$, then the have the same 
  Amitsur subgroups inside $\Br(k)$ 
  and in particular, the same period. \qed
\end{Theorem}

\begin{Remark}
  \label{rem: the amitsur conjecture}
  In general, it is not true that two Brauer--Severi varieties of the same dimension
  and the same Amitsur subgroup are isomorphic.
  We refer to Remark \ref{rem: amitsur remark} for an example arising from 
  a Cremona transformation of Brauer--Severi surfaces.
  However, Amitsur asked whether 
  two Brauer--Severi varieties of the same dimension with the same Amitsur
  subgroup are birationally equivalent.
\end{Remark}

In our applications to del~Pezzo surfaces below,
we will only need the following easy and probably well-known corollary.

\begin{Corollary}
  \label{cor: zero cycle}
  Let $P$ be a Brauer--Severi variety over $k$.
  If there exists a zero-cycle on $P$, whose
  degree is prime to $(\dim(P)+1)$,
  then $P$ is  is trivial.
\end{Corollary}

\prf
Since ${\rm Am}(P)\iso\ZZ/{\rm per}(P)\ZZ$ and its order 
divides $(\dim(P)+1)$,
Lemma \ref{lem: order amitsur} and the assumptions imply
${\rm Am}(P)=0$.
Thus, ${\rm per}(P)=1$, and then, $P$ is trivial.
\qed\medskip

We end this section by mentioning another important invariant of a Brauer--Severi variety 
$P$ over $k$, namely, its {\em index}, denoted ${\rm ind}(P)$.
We refer to \cite[Chapter 4.5]{Gille Szamuely} for the precise definition and note
that it is equal to the smallest degree of a finite separable field extension $K/k$ 
such that $P_K$ is trivial, as well as to the greatest common divisor 
of the degrees of all finite separable
field extensions $K/k$ such that $P_K$ is trivial.
By a theorem of Brauer, the period divides the index, and they have the same 
prime factors, see \cite[Proposition 4.5.13]{Gille Szamuely}.

\section{Morphisms to Brauer--Severi varieties}

This section contains Theorem \ref{thm: main}, the main observation of this article
that describes morphisms from a proper variety $X$ over a field $k$ to 
Brauer--Severi varieties in terms of classes in of $\Pic_{(X/k)({\rm fppf})}(k)$.
We start by extending classical notions for invertible sheaves to such
classes, and then, use these notions to phrase and prove Theorem \ref{thm: main}.
As immediate corollaries, we obtain two classical results of Kang and Ch\^{a}telet
on the geometry of Brauer--Severi varieties.

\subsection{Splitting fields, globally generated and ample classes}
Before coming to the main result of this section, we introduce the following.

\begin{Definition}
 \label{def: globally generated}
 Let $X$ be a proper variety over $k$ and ${\cal L}\in\Pic_{(X/k)({\rm fppf})}(k)$.
 \begin{enumerate}
  \item A {\em splitting field} for $\cal L$ is a field extension $K/k$ such that
     ${\cal L}\otimes_kK$ lies in $\Pic(X_K)$, i.e., arises from an invertible sheaf 
     on $X_K$.
   \item The class $\cal L$ is called {\em globally generated} (resp. {\em ample}, resp. {\em very ample}) 
     if there exists a splitting
     field $K$ for ${\cal L}$ such that ${\cal L}\otimes_kK$ is globally generated (resp. ample, resp. very ample)
     as an invertible sheaf on $X_K$.
 \end{enumerate}
\end{Definition}

From the short exact sequence in Proposition \ref{prop: delta}, 
it follows that if $K$ is a splitting field for the class ${\cal L}$, then there exists precisely one invertible sheaf 
on $X_K$ up to isomorphism  that corresponds to this class.
The following lemma shows that these notions are independent of the choice
of a splitting field of the class $\cal L$.

\begin{Lemma}
  \label{well defined lemma}
  Let $X$ be a proper variety over $k$ and ${\cal L}\in\Pic_{(X/k)({\rm fppf})}(k)$.
  \begin{enumerate}
   \item There exists a splitting field for $\cal L$ that is a finite and separable extension $k$, and it can also
     chosen to be Galois over $k$.
   \item Let $K$ and $K'$ be splitting fields for ${\cal L}$.
     Then ${\cal L}\otimes_k K\in\Pic(X_K)$ is globally generated (resp. ample, resp. very ample)
     if and only if  ${\cal L}\otimes_k K'\in\Pic(X_{K'})$ is globally generated (resp. ample, resp. very ample).
  \end{enumerate}
\end{Lemma}

\prf
To simplify notation in this proof, we set ${\cal L}_K:={\cal L}\otimes_k K$.

Let $K$ be a finite and separable extension of $k$, such that $\delta({\cal L})\in\Br(k)$
lies in $\Br(K|k)$, where $\delta$ is as in Proposition \ref{prop: delta}.
Then, $\delta({\cal L}_K)=0$, i.e., ${\cal L}_K$ 
comes from an invertible sheaf on $X_K$.
In particular, $K$ is a splitting field for $\cal L$, which is a finite and separable extension of $k$.
Passing to the Galois closure of $K/k$, we obtain a splitting field for $\cal L$ that is a finite
Galois extension of $k$.
This establishes claim (1).

Claim (2) is a well-known application of flat base change, 
but let us recall the arguments for the reader's convenience:
By choosing a field extension of $k$ that contains both $K$ and $K'$, we reduce to the case
$k\subseteq K\subseteq K'$.
We have $H^0(X_K,{\cal L}_K)\otimes_K K'\iso H^0(X_{K'},{\cal L}_{K'})$ by flat base change
for cohomology, from which it is easy to see that
${\cal L}_K$ is globally generated if and only if ${\cal L}_{K'}$ is so.
Next, if ${\cal L}_K$ is very ample, then its global sections give rise to
a closed immersion $X_K\to\PP^n_K$ for some $n$.
After base change to $K'$, we obtain a closed embedding $X_{K'}\to\PP^n_{K'}$ which
corresponds to the  global sections of ${\cal L}_{K'}$, and so,
also ${\cal L}_{K'}$ is very ample.
Conversely, if ${\cal L}_{K'}$ is very ample, then it is globally generated, and thus,
${\cal L}_K$ is globally generated by what we just established,
and thus, gives rise to a morphism $\varphi_K:X_K\to\PP^n_K$.
By assumption and flat base change, $\varphi_{K'}$ is a closed embedding, and
thus, $\varphi_K$ is a closed embedding, and ${\cal L}_K$ is very ample.
From this, it also follows that ${\cal L}_K$ is ample if and only if ${\cal L}_{K'}$ is.
\qed\medskip

\begin{Remark}
  \label{rem: explain delta} 
  Let $X$ be a proper variety over $k$ and let
  $$
    \begin{array}{ccccc}
     \delta &:& \Pic_{(X/k)({\rm fppf})}(k) &{\longrightarrow}& \Br(k)
    \end{array}
  $$
  be as in Proposition \ref{prop: delta}.
 We are now in a position to describe $\delta$ explicitly.
  \begin{enumerate}
   \item First, and more abstractly: given a class ${\cal L}\in\Pic_{(X/k)({\rm fppf})}(k)$, we can choose 
     a splitting field $K$ that is a finite extension $k$. 
     Thus, $\Spec K\to \Spec k$ is an fppf cover, 
     the class ${\cal L}\otimes_k K$ comes with an fppf descent datum, and it
     arises from an invertible sheaf ${\cal M}$ on $X_K$.
     The crucial point is that the descent datum is for a class in $\Pic(X_K)$, where isomorphism classes of
     invertible sheaves are identified.
     In order to turn this into a descent datum for the invertible sheaf $\cal M$, we have to choose 
     isomorphisms,
     which are only unique up to a $\GG_m={\rm Aut}({\cal M})$-action, and we obtain
     a $\GG_m$-gerbe that is of class $\delta({\cal L})\in H^2_{{\rm fppf}}(\Spec k,\GG_m)=\Br(k)$.
     This gerbe is neutral if and only if $\delta({\cal L})=0$.
     This is equivalent to being able to extend the descent datum for the class ${\cal L}\otimes_k K$
     to a descent datum for the invertible sheaf ${\cal M}$.
    \item Second, and more concretely: given a class ${\cal L}\in\Pic_{(X/k)({\rm fppf})}(k)$, 
     we can choose a splitting field $K$ that is a finite Galois extension of $k$, say with Galois group $G$.
     Thus, the class ${\calL}\otimes_k K$ arises from an invertible 
     sheaf $\cal M$ on $X_K$ and lies in $\Pic_X(K)^G$
     and we can choose isomorphisms
     $$
      \imath_g \,:\, g^\ast {\cal M} \,\stackrel{\iso}{\longrightarrow}\, {\cal M},
     $$
     which are unique up to a $\GG_m$-action.
     In particular, they may fail to form a Galois descent datum for $\cal M$,
     and the failure of turning $\{\imath_g\}_{g\in G}$ into a Galois descent datum for $\cal M$
     gives rise to a cohomology class $\delta({\cal L})\in \Het{2}(\Spec k,\GG_m)=\Br(k)$.
     More precisely, this class lies in the subgroup $\Br(K|k)$ of $\Br(k)$.
 \end{enumerate}
\end{Remark}

The following is an analog for Brauer--Severi varieties of the classical correspondence
between morphisms to projective space and globally generated invertible sheaves as explained, 
for example, in \cite[Theorem II.7.1]{Hartshorne}, see also Remark \ref{rem: trivial case} 
below.

\begin{Theorem}
  \label{thm: main}
   Let $X$ be a proper variety over a field $k$.
   \begin{enumerate}
    \item Let $\varphi:X\to P$ be a morphism to a Brauer--Severi variety $P$ over $k$,
      and consider the induced homomorphism of abelian groups
      $$
        \begin{array}{ccccc}
          \varphi^* &:& \Pic_{(P/k)({\rm fppf})}(k) &\to& \Pic_{(X/k)({\rm fppf})}(k).
        \end{array}
     $$
      Then, ${\cal L}:=\varphi^*\OO_P(1)$ is a globally generated class with
      $$
        \begin{array}{ccc}
          \delta({\cal L}) \,=\, [P] &\in& \Br(k),
        \end{array}
       $$
      where $\delta$ is as in Proposition \ref{prop: delta}.
      If $\varphi$ is a closed immersion, then $\cal L$ is very ample.
    \item Let ${\cal L}\in\Pic_{(X/k)({\rm fppf})}(k)$ be a globally generated class.
      If $K$ is a splitting field, then the morphism to projective space over $K$
      associated to the complete linear system 
      $|{\cal L}\otimes_kK|$ descends to morphism over $k$
      $$
        \begin{array}{ccccc}
          |{\cal L}| &:& X &\to& P,
        \end{array}
     $$
     where $P$ is a Brauer--Severi variety over $k$ with
     $\delta({\cal L})=[P]$.
     If ${\cal L}$ is very ample, then $|{\cal L}|$ is a closed immersion.
  \end{enumerate}
\end{Theorem}

\prf
Let $\varphi:X\to P$ and $\cal L$ be as in (1).
Then, we have $\delta({\cal L})=\delta(\OO_P(1))=[P]\in\Br(k)$,
where the first equality follows from functoriality of the exact
sequence in Proposition \ref{prop: delta}, and the second from 
Theorem \ref{thm: brauer severi picard}.
Let $K$ be a splitting field for $\cal L$, and let  $\cal M$ be 
the invertible sheaf corresponding to ${\cal L}\otimes_k K$ on $X_K$.
Being an invertible sheaf, we have
$\delta({\cal M})=0\in\Br(K)$, which implies that
the morphism $\varphi_K:X_K\to P_K$ maps to a Brauer--Severi variety
of class $[P_K]=\delta({\cal M})=0$, i.e., $P_K\iso\PP_K^n$.
By definition and base change, we obtain
${\cal M}\iso\varphi_K^*(\OO_{\PP_K^n}(1))$.
Thus, $\cal M$ is globally generated (as an invertible sheaf), which 
implies that ${\cal L}\in\Pic_{(X/k)({\rm fppf})}(k)$ is globally generated in 
the sense of Definition \ref{def: globally generated}.
Moreover, if $\varphi$ is a closed immersion, then so is $\varphi_K$,
which implies that ${\cal M}\in\Pic(X_K)$ is very ample (as an invertible sheaf), 
and thus, ${\cal L}\in\Pic_{(X/k)({\rm fppf})}(k)$ is very ample in the sense of 
Definition \ref{def: globally generated}.
This establishes claim (1)

To establish claim (2), let ${\cal L}\in\Pic_{(X/k)({\rm fppf)}}(k)$ be globally generated.
By Lemma \ref{well defined lemma}, there exists a splitting field $K'$ for ${\cal L}$ that is
a finite Galois extension of $k$, say with Galois group $G$.
Thus, ${\cal L}\otimes_kK'$ corresponds to an invertible sheaf ${\cal M}$ on $X_{K'}$, 
whose isomorphism class lies in $\Pic_{X}(K')^G$,
see Proposition \ref{prop: picard in geometry}.

If $f:X\to\Spec k$ is the structure morphism, then $(f_{K'})_\ast{\cal M}$ is
a finite-dimensional $K'$-vector space.
By our assumptions on global generation we obtain a morphism over $K'$
$$
  |{\cal M}| \,:\,X_{K'} \,\to\, \PP((f_{K'})_\ast{\cal M}).
$$
As explained in Remark \ref{rem: explain delta}.(2), there exist isomorphisms
$\{\imath_g:g^*{\cal M}\to{\cal M}\}_{g\in G}$ that are unique up to a $\GG_m$-action.
In particular, we obtain a well-defined $G$-action on 
$\PP((f_{K'})_\ast{\cal M})$, and the morphism defined
by $|{\cal M}|$ is $G$-equivariant.
Taking the quotient by $G$, we obtain a morphism over $k$
$$
   |{\cal L}|\,:\,X\to P.
$$
Since $P_K$ is isomorphic to $\PP((f_{K'})_\ast{\cal M})$, we see that $P$
is a Brauer--Severi variety over $k$ and, as observed by Grothendieck in \cite[Section (5.4)]{Grothendieck Brauer}, 
we have $\delta({\cal L})=[P]$ in $\Br(k)$. 

Finally, let $K$ be an arbitrary splitting field for ${\cal L}$.
Let $\varphi:X\to P$ be the previously constructed morphism and choose an
extension field $\Omega$ of $k$ that contains $K$ and $K'$.
Then, ${\cal L}\otimes_k\Omega$ is an invertible sheaf on $X_\Omega$, globally generated by 
Lemma \ref{well defined lemma}, and, since $k\subseteq K'\subseteq\Omega$, the morphism
associated to $|{\cal L}\otimes_k\Omega|$ is % by construction of $\varphi_{K'}$ 
equal to $\varphi_\Omega=(\varphi_{K'})_\Omega:X_\Omega\to P_\Omega$.
Since $K$ is a splitting field for ${\cal L}$, it is also a splitting field for $P_{K}$ (see the argument in the 
proof of claim (1)), and in particular, $P_{K'}$ is a trivial Brauer--Severi variety.
We have ${\cal L}\otimes_k\Omega\iso \varphi_\Omega^*\OO_{P_\Omega}(1)$, from which 
we deduce ${\cal L}\otimes_k K\iso \varphi_{K}^*\OO_{P_{K}}(1)$, as well as that
$\varphi_{K}$ is the morphism associated to $|{\cal L}\otimes_k K|$.
In particular, the morphism associated to 
$|{\cal L}\otimes_k K|$ descends to $\varphi:X\to P$, where $P$ is
a Brauer--Severi variety of class $\delta({\cal L})$.
This establishes claim (2).
\qed\medskip

\begin{Remark}
 \label{rem: trivial case}
 Let us note the following.
 \begin{enumerate}
  \item The construction of a Brauer--Severi variety over $k$ 
    from a globally generated class in $\Pic_{(X/k)({\rm fppf})}(k)$ (in our terminology) 
    is due to Grothendieck in \cite[Section (5.4)]{Grothendieck Brauer}.
  \item In Theorem \ref{thm: main}.(2), we only considered complete linear systems.
    We leave it to the reader to show the following generalization: 
    Given a class ${\cal L}\in\Pic_{(X/k)({\rm fppf})}(k)$,
    a splitting field $K$ that is finite and Galois over $k$ with Galois group $G$,
    and $V\subseteq H^0(X_K,{\cal L}\otimes_kK)$ a $G$-stable $K$-linear subspace, 
    whose global sections generate ${\cal L}\otimes_kK$, we can descend the
    morphism $X_K\to\PP(V)$ to a morphism $X\to P'$, where $P'$ is a Brauer--Severi
    variety over $k$ of class $[P']=\delta({\cal L})\in\Br(k)$.
  \item If $X$ in Theorem \ref{thm: main} has a $k$-rational point, i.e., 
    $X(k)\neq\emptyset$, then we recover the well-known correspondence between
    morphisms to projective space and globally generated invertible sheaves:
    \begin{enumerate}
     \item Then, $\delta\equiv0$ and every class in $\Pic_{(X/k)({\rm fppf})}(k)$
         comes from an invertible sheaf on $X$ by Proposition \ref{prop: delta},
     \item and since every morphism $\varphi:X\to P$ 
         gives rise to a $k$-rational point on $P$, i.e., $P$ is a trivial
        Brauer--Severi variety.
   \end{enumerate}
 \end{enumerate}
 \end{Remark}

\subsection{Two classical results on Brauer--Severi varieties}
As our first corollary and application, we recover the following theorem of Kang \cite{Kang}, 
see also \cite[Theorem 5.2.2]{Gille Szamuely}, which is a Brauer--Severi variety
analog of Veronese embeddings of projective spaces.

\begin{Corollary}[Kang]
 \label{cor: kang}   
 Let $P$ be a Brauer--Severi variety of period ${\rm per}(P)$ over $k$.
 Then, the class of $\OO_{P}({\rm per}(P))$ arises from 
 a very ample invertible sheaf on $P$ and gives rise to an embedding 
 $$
   |\OO_P({\rm per}(P))| \,:\, P \,\to\, \PP^N_k,
    \mbox{ \quad where \quad } N\,=\, \binom{\dim(P)+{\rm per}(P)}{{\rm per}(P)} .
  $$
  After base change to a splitting field $K$ of $P$, this embedding becomes
  the ${\rm per}(P)$-uple Veronese embedding of $\PP^{\dim(P)}_K$ into
  $\PP^N_K$.
\end{Corollary}

\proof
If  $n\geq1$, then $\OO_P(n)$ is very ample 
in the sense of Definition \ref{def: globally generated}, and thus,
defines an embedding into a Brauer--Severi variety $P'$ over $k$.
Over a splitting field of $P$, this embedding becomes the $n$-uple
Veronese embedding.
Since $\delta(\OO_P(1))=[P]\in\Br(k)$ and this element of order
${\rm per}(P)$, we see that if ${\rm per}(P)$ divides $n$, then
$\OO_P(n)$ is an invertible sheaf on $P$ and
$P'$ is a trivial Brauer--Severi variety.
\qed\medskip

\begin{Example}
 \label{example: BS curve}
 Let $X$ be a smooth and proper variety of dimension one over $k$.
 If $\omega_X^{-1}$ is ample, then it is a curve of genus
 $g(X)=h^0(X,\omega_X)=0$.
 Thus, $X$ is isomorphic to $\PP^1$ over $\overline{k}$, i.e., 
 $X$ is a Brauer--Severi curve.
There exists a unique class $\calL\in\Pic_{(X/k)({\rm fppf})}(k)$
 with $\calL^{\otimes2}\iso\omega_X^{-1}$, and it gives rise to
 an isomorphism $|\calL|:X\to P$, where $P$ is a Brauer--Severi curve with
 $\delta(\calL)=[P]\in\Br(k)$.
 Moreover, $\calL^{\otimes2}\iso\omega_X^{-1}$ is an invertible sheaf on $X$ 
 that defines an embedding $|\omega_X^{-1}|:X\to\PP^2_k$ as a plane conic.
\end{Example}

A subvariety $X\subseteq P$ of a Brauer--Severi
variety $P$ over $k$ is called {\em twisted linear} if $X_{\overline{k}}$
is a linear subspace of $P_{\overline{k}}$. 
As second application, we recover the following theorem of Ch\^{a}telet, 
see \cite[Section 5.3]{Gille Szamuely},  and it follows from a Brauer--Severi variety
analog of Segre embeddings of products of projective spaces.

\begin{Corollary}[Ch\^{a}telet]
  \label{cor: chatelet}
  Let $P_1$ and $P_2$ be two Brauer--Severi varieties over $k$
  of dimension $d_1$ and $d_2$, respectively.
  \begin{enumerate}
   \item If $P_1$ is a twisted linear subvariety of $P_2$, then
    $[P_1]=[P_2]\in\Br(k)$.
    \item If $[P_1]=[P_2]\in\Br(k)$, then there exists a 
     Brauer--Severi variety $P$ over $k$, such that $P_1$ and $P_2$ 
     can be embedded as twisted-linear subvarieties into $P$.
  \end{enumerate}
\end{Corollary} 

\prf
If $\varphi:P_1\into P_2$ is a twisted-linear subvariety,
then $\varphi^*\OO_{P_2}(1)=\OO_{P_1}(1)\in\Pic_{(P_1/k)({\rm fppf})}(k)$.
We find $[P_1]=\delta(\OO_{P_1}(1))=\delta(\OO_{P_2}(1))=[P_2]$
by functoriality of the exact sequence of Proposition \ref{prop: delta},
and (1) follows.

Next, we show (2).
By Theorem \ref{thm: main}, there exists an embedding $\varphi$ 
of $P_1\times\PP_k^{d_2}$ into a Brauer--Severi variety $P$ 
of dimension $N:=(d_1+1)(d_2+1)-1=d_1d_2+d_1+d_2$ over $k$ associated to the
class $\OO_{P_1}(1)\boxtimes\OO_{\PP_k^{d_2}}(1)$.
Over a splitting field of $P_1$, this embedding becomes the Segre embedding
of $\PP^{d_1}\times\PP^{d_2}$ into $\PP^N$.
If $x$ is a $k$-rational point of $\PP_k^{d_2}$, then 
$\varphi(P_1\times\{x\})$ realizes $P_1$ as twisted-linear subvariety of $P$
and we have $[P]=[P_1]\in\Br(k)$ by claim (1).
Similarly, we obtain an embedding of $P_2$ as twisted-linear
subvariety into a Brauer--Severi variety $P'$ of dimension 
$N$ over $k$ of class $[P']=[P_2]\in\Br(k)$.
Since $[P]=[P']\in\Br(k)$ and $\dim(P)=\dim(P')$, 
we find $P\iso P'$ by Corollary \ref{cor: isomorphic BS}
and (2) follows.
\qed\medskip

\section{Del~Pezzo surfaces}

For the remainder of this article, we study del~Pezzo surfaces with a view towards 
Brauer--Severi varieties.
Most, if not all, results of these sections are known in some form or another to
the experts. 
However, our more geometric approach, as well as some of the proofs, are new.

Let us first recall some classical results about del~Pezzo surfaces, and
refer the reader to \cite[Chapter IV]{Manin} or the surveys
\cite{CT survey}, \cite{Varilly}, \cite{Poonen} for details, proofs, and references.
For more results about the classification of geometrically rational 
surfaces, see \cite{Manin surfaces} and \cite{Iskovskih}.

\begin{Definition}
  A {\em del~Pezzo surface} is a smooth and proper variety $X$ of dimension two
  over a field $k$ such that $\omega_X^{-1}$ is ample.
  The {\em degree} of a del~Pezzo surface is the self-intersection number of $\omega_X$.
\end{Definition}

In arbitrary dimension, smooth and proper varieties $X$ over $k$ with
ample $\omega_X^{-1}$ are called {\em Fano varieties}.
As discussed in Example  \ref{example: BS curve}, Fano varieties of dimension one over $k$
are the same as Brauer--Severi curves over $k$.

\subsection{Geometry}
The degree $d$ of a del~Pezzo surface $X$ over a field $k$ 
satisfies $1\leq d\leq 9$.
Set $\overline{X}:=X_{\overline{k}}$.
We will say that $X$ is {\em of product type} if
$$
 \begin{array}{ccc}
  \overline{X}  &\iso& \PP^1_{\overline{k}}\times\PP^1_{\overline{k}},
 \end{array}
$$
in which case we have $d=8$.
If $X$ is not of product type, then there exists a birational morphism
$$
 \begin{array}{ccccc}
   \overline{f} &:&\overline{X}&\to&\PP^2_{\overline{k}}
 \end{array}
$$
that is a blow-up of $(9-d)$ closed points $P_1,...,P_{9-d}$ in general position, i.e.,
no $3$ of them lie on a line, no $6$ of them lie on a conic, and there is no cubic
through all these points having a double point in one of them.
In particular, if $d=9$, then $\overline{f}$ is an isomorphism and
$X$ is a Brauer--Severi surface over $k$.

\subsection{Arithmetic}
\label{subsec: dP arithmetic}
By the previous discussion and Lemma \ref{lem: picard rank}, the
{\em N\'eron--Severi rank} of a del~Pezzo surface $X$ of degree $d$ 
over $k$ satisfies
$$
    1\,\leq\,\rho(X)\,:=\, {\rm rank}\,\Pic(X)\,=\,{\rm rank}\,\Pic_{(X/k)({\rm fppf})}(k) \,\leq\,10-d,
$$
and $\rho(X_{\overline{k}})=10-d$.

The following result about geometrically rational surfaces 
allows using methods from Galois theory even if the ground field $k$ is 
not perfect.
This result is particularly useful in proofs, see also the 
discussion in \cite[Section 1.4]{Varilly}.
In particular, it applies to del~Pezzo surfaces.

\begin{Theorem}[Coombes$+\varepsilon$]
 \label{thm: Coombes}
  Let $X$ be a smooth and proper variety over $k$ 
  such that $X_{\overline{k}}$ is birational to $\PP^2_{\overline{k}}$.
  Then,
  \begin{enumerate}
    \item $X_{k^{\rm sep}}$ is birationally equivalent to $\PP^2_{k^{\rm sep}}$ via a sequence
      of blow-ups in points in $k^{\rm sep}$-rational points and their inverses.
     \item The natural map $\Pic_X(k^{\rm sep})\to\Pic_X({\overline{k}})$ is an
      isomorphism.
  \end{enumerate}
\end{Theorem}

\prf
Assertion (1) is the main result of \cite{Coombes}.
Clearly, assertion (2) holds for projective space over any field.
Next, let $Y$ be a variety that is smooth and proper over $k^{\rm sep}$,
$\widetilde{Y}\to Y$ be the blow-up of a $k^{\rm sep}$-rational point, and let
$E\subset\widetilde{Y}$ be the exceptional divisor.
Then, $\Pic_{\widetilde{Y}}(K)=\Pic_Y(K)\oplus\ZZ\cdot E$ for $K=k^{\rm sep}$,
as well as for $K=\overline{k}$.
Using (1) and these two observations, assertion (2) follows.
\qed\medskip

We will also need the following useful observation, due to Lang \cite{Lang} and
Nishimura \cite{Nishimura}, which implies that having a $k$-rational point
is a birational invariant of smooth and proper varieties over $k$.
We refer to \cite[Section 1.2]{Varilly} for details and proof.

\begin{Lemma}[Lang--Nishimura] 
  \label{lem: Lang}
  Let $X\dashrightarrow Y$ be a rational map of varieties over $k$, such that
  $X$ is smooth over $k$, and such that $Y$ is proper over $k$.
  If $X$ has a $k$-rational point, then so has $Y$.\qed
\end{Lemma}

Moreover, we have already seen that a Brauer--Severi variety $P$ over $k$
is isomorphic to projective space over $k$
if and only if $P$ has a $k$-rational point, and we refer the interested reader to \cite{BS algorithm} 
for an algorithm to decide whether a Brauer--Severi surface has a $k$-rational point.
In Definition \ref{def: amitsur group}, we defined the Amitsur group and showed its
birational invariance in Proposition \ref{prop: amitsur is birational invariant}.
Using Iskovskih's classification \cite{Iskovskih} of 
geometrically rational surfaces, we obtain the following list
and refer to \cite[Proposition 5.2]{CTKMdP6} for details and proof.

\begin{Theorem}[Colliot-Th\'el\`ene--Karpenko--Merkurjev]
  Let $X$ be a smooth and proper variety over a perfect field $k$ such that
   $X_{\overline{k}}$ is birationally equivalent to $\PP^2_{\overline{k}}$.
  Then, ${\rm Am}(X)$ is one of the following groups
  $$
   0,\mbox{ \quad }\ZZ/2\ZZ,\mbox{ \quad } (\ZZ/2\ZZ)^2,\mbox{ \quad and \quad }\ZZ/3\ZZ.
  $$
\end{Theorem}

We will see explicit examples of all these groups arising as Amitsur groups of
del~Pezzo surfaces in the next sections.

We now introduce another important invariant.
Namely, if $G_k$ denotes the absolute Galois group of $k$, 
and $H\subseteq G_k$ is a closed subgroup,
then we consider for a smooth and projective variety $X$ over $k$ 
the group cohomology 
$$
   H^1\left(H,\, \Pic_{X/k}({k^{\rm sep}})\right),
$$
which is an abelian torsion group.
If $b_1(X)=0$, then $\Pic_{X/k}(k^{\rm sep})$ is finitely generated 
by Lemma \ref{lem: picard rank}
and then, $H^1(H,\, \Pic_{X/k}({k^{\rm sep}}))$ is a finite abelian group.
Moreover, if $X_{k^{\rm sep}}$ is a rational surface, 
then $\Br'(X_{k^{\rm sep}})=0$ (see, for example, \cite[Theorem 42.8]{Manin}
or \cite{Milne Brauer}) and an appropriate 
Hochschild--Serre spectral sequence yields an exact sequence
$$
  0\,\to\,\Br'(X)/\Br(k)\,\stackrel{\alpha}{\longrightarrow}\,
  H^1\left(G_k,\, \Pic_{X/k}({k^{\rm sep}})\right)
  \,\to\, H^3(G_k,(k^{\rm sep})^\times).
$$
Moreover, if $k$ is a global field, then 
the term on the right is zero by a theorem of Tate
(see, for example, \cite[Chapter VIII.3]{Neukirch}), thus, 
$\alpha$ is an isomorphism, 
and we obtain an interpretation of this cohomology group
in terms of Brauer groups, see \cite[Section 3.4]{Varilly}.

\begin{Lemma}
  \label{lem: h1 brauer-severi}  
  If $P$ is a Brauer--Severi variety over $k$, then
  $$
     H^1\left(H,\, \Pic_{P/k}({k^{\rm sep}})\right)\,=\,0
  $$
  for all closed subgroups $H\subseteq G_k$.
\end{Lemma}

\prf
Since $\Pic_{P/k}({k^{\rm sep}})\iso\ZZ\cdot \OO_P(1)$ and since
$G_k$ acts trivially on the class $\OO_P(1)$, the desired $H^1$
is isomorphic to ${\rm Hom}(H,\ZZ)$,
see \cite[Chapter III.1, Exercise 2]{Brown}, for example.
This is zero since $H$ is a profinite group and the homomorphisms
to $\ZZ$ are required to be continuous.
\qed\medskip

In Proposition \ref{prop: amitsur is birational invariant},
we established birational invariance of ${\rm Am}(X)$.
The following result of Manin \cite[Section 1 of the Appendix]{Manin} shows 
that also the above group cohomology groups are a birational invariants.

\begin{Theorem}[Manin]
  \label{thm: birational invariance of h1}
  For every closed subgroup $H\subseteq G_k$, the group
  $$ 
     H^1\left(H,\, \Pic_{X/k}({k^{\rm sep}})\right) 
  $$
  is a birational invariant of smooth and projective varieties over $k$.
  \qed
\end{Theorem}

\begin{Remark}
  Every birational map between smooth and projective surfaces
  can be factored into a sequence of blow-ups in closed points,
  see \cite[Chapter III]{Manin}.
  Using this, one can give very explicit proofs of
  Proposition \ref{prop: amitsur is birational invariant} and 
  Theorem \ref{thm: birational invariance of h1} 
  in dimension $2$. 
  (For such a proof of Theorem \ref{thm: birational invariance of h1} 
  in dimension $2$, see the proof of \cite[Theorem 29.1]{Manin}.)
\end{Remark}

\subsection{Hasse principle and weak approximation}
\label{subsec: Hasse}
For a global field $K$, i.e., a finite extension of $\QQ$ or of $\FF_p(t)$, we denote by
$\Omega_K$ the set of its places, including the infinite ones if $K$ is of 
characteristic zero.
A class $\calC$ of varieties over $K$ satisfies

\begin{enumerate}
 \item the {\em Hasse principle}, if for every $X\in\calC$ we have $X(K)\neq\emptyset$ if and only if 
    $X(K_\nu)\neq\emptyset$ for all $\nu\in\Omega_K$. 
    Moreover, $\calC$ satisfies
 \item {\em weak approximation}, if the diagonal embedding
    $$
        X(K) \to \prod_{\nu\in\Omega_K} X(K_\nu)
    $$
    is dense for the product of the $\nu$-adic topologies.
\end{enumerate}
If $\calC$ satisfies weak approximation, then it obviously also satisfies the Hasse principle, 
but the converse need not hold.
For example, Brauer--Severi varieties over $K$ satisfy 
the Hasse principle by a theorem of Ch\^{a}telet \cite{Chatelet}, as well as
weak approximation.
However, both properties may fail for del~Pezzo surfaces over $K$, and
we refer to \cite{Varilly} for an introduction to this topic.
We end this section by noting that the obstruction to 
a class $\Pic_{(X/K)({\rm fppf})}(K)$ coming from $\Pic_X(K)$
satisfies the Hasse principle.

\begin{Lemma}
  Let $X$ a proper variety over a global field $K$ and let
  $\calL\in\Pic_{(X/K)({\rm fppf})}(K)$.
  Then, the following are equivalent
  \begin{enumerate}
  \item $0=\delta({\cal L})\in\Br(K)$ , and
  \item $0=\delta(\calL\otimes_K K_\nu)\in\Br(K_\nu)\mbox{ for all }\nu\in\Omega_K$.
  \end{enumerate}
\end{Lemma}

\prf
A class in $\Br(K)$ is zero if and only if its image in $\Br(K_\nu)$ is zero
for all $\nu\in\Omega_K$ by the Hasse principle for the Brauer group.
From this, and functoriality of the exact sequence from Proposition \ref{prop: delta},
the assertion follows.
\qed\medskip

For example, if $X(K_\nu)\neq\emptyset$ for all $\nu\in\Omega_X$, then $\delta$ is 
the zero map by Proposition \ref{prop: delta} and this lemma. 
In this case, every class in $\Pic_{(X/K)({\rm fppf})}(K)$
comes from an invertible sheaf on $X$.

\section{Del~Pezzo surfaces of product type}

In this section, we classify degree $8$ del~Pezzo surfaces of product type over $k$, 
i.e., surfaces $X$ over $k$ with $X_{\overline{k}}\iso\PP^1_{\overline{k}}\times\PP^1_{\overline{k}}$,
in terms of Brauer--Severi varieties.

First, for $\PP^1_k\times\PP^1_k$, the anti-canonical embedding  can be written as composition
of Veronese- and Segre-maps as follows
$$
 \begin{array}{ccccccc}
 |-K_{\PP^1_k\times\PP^1_k}| &:& \PP^1_k\times\PP^1_k 
 &\stackrel{\nu_2\times\nu_2}{\longrightarrow}& \PP^2_k\times\PP^2_k &\stackrel{\sigma}{\longrightarrow}& \PP^8_k\,.
 \end{array}
$$
Next, the invertible sheaf $\omega_{\PP^1_k\times\PP^1_k}^{-1}$ 
is uniquely $2$-divisible in the Picard group, and we obtain an embedding as a smooth quadric
$$
 \begin{array}{ccccc}
 |{\scriptstyle -\frac{1}{2}}K_{\PP^1_k\times\PP^1_k}| &:&
 \PP^1_k\times\PP^1_k &\stackrel{\sigma}{\longrightarrow}& \PP^3_k\,.
 \end{array}
$$
Now, let $X$ be a degree $8$ del~Pezzo surface of product type over $k$.
Then, the anti-canonical linear system  yields an embedding of $X$ as a surface
of degree $8$ into $\PP^8_k$.
However, the ``half-anti-canonical linear system'' exists in general only as a morphism 
to a Brauer--Severi threefold as the following result shows.

\begin{Theorem}
 \label{thm: product type}  
 Let $X$ be a degree $8$ del~Pezzo surface of product type over a field $k$.
 Then, there exist a unique class ${\cal L}\in\Pic_{(X/k)({\rm fppf})}(k)$ with
 ${\cal L}^{\otimes 2}\iso\omega_X^{-1}$ and an embedding
  $$\begin{array}{ccccc}
         |\calL| &:& X &\into&P
      \end{array}
  $$
  into a Brauer--Severi threefold $P$ over $k$ with Brauer class
  $$\begin{array}{ccc}
         \delta({\cal L}) \,=\, [P] &\in&\Br(k),
       \end{array}
  $$
  and such that $X_{\overline{k}}$ is a smooth quadric in $P_{\overline{k}}\iso\PP^3_{\overline{k}}$.  
  Moreover, $X$ is rational if and only if $X$ has a $k$-rational point.
  In this case, we have $P\iso\PP^3_k$.
\end{Theorem}

\prf
To simplify notation, set $L:=k^{\rm sep}$.
We have $X(L)\neq\emptyset$, for example, by \cite[Proposition A.1.1]{Gille Szamuely},
as well as  $\Pic(X_L)\iso\Pic(X_{\overline{k}})\iso\ZZ^2$ by 
Theorem \ref{thm: Coombes}.
The classes $(1,0)$ and $(0,1)$ of $\Pic(X_L)$ give rise to two morphisms
$X_L\to\PP^1_L$, and we obtain an isomorphism 
$X_L\iso\PP^1_L\times\PP^1_L$.
By abuse of notation, we re-define $\overline{X}$ to be $X_L$.
Next, the absolute Galois group $G_k$ acts trivially on the
canonical class $(-2,-2)$, and thus, the $G_k$-action on  $\ZZ(1,1)\subset\ZZ^2$ is trivial.
By Proposition \ref{prop: picard in geometry}, we have
$\Pic_{X/k}(K)^{G_k}\iso\Pic_{(X/k)({\rm fppf})}(k)$, and,
since $(1,1)\in\ZZ^2$ is $G_k$-invariant, 
the unique invertible sheaf $\calL$ on $\overline{X}$ with
$\calL^{\otimes 2}\iso\omega_{\overline{X}}^{-1}$ descends to a class in
$\Pic_{(X/k)({\rm fppf})}(k)$.
Over $L$, the class $\calL$ is very ample and defines an embedding of 
$\overline{X}$ as smooth quadric surface into $\PP^3_L$.
Thus, by Theorem \ref{thm: main}, we obtain an embedding 
$|\calL|:X\into P$, where $P$ is a Brauer--Severi threefold over $k$ 
with $\delta(\calL)=[P]\in\Br(k)$.

Finally, if $X$ is rational, then it has a $k$-rational point, and then,
also $P$ has a $k$-rational point, i.e., $P\iso\PP^3_k$.
Conversely, if there exists a $k$-rational point $x\in X$, then $X$ is a quadric in $\PP^3_k$,
and projection away from $x$ induces a birational map
$X\dashrightarrow\PP^2_k$.
\qed\medskip

Next, we establish an explicit classification of degree $8$ 
del~Pezzo surfaces of product type in terms of the N\'eron--Severi
rank $\rho$ and Brauer--Severi curves.
To simplify notation in the sequel, let us recall the definition
of contracted products.
If a finite group $G$ acts on a scheme $X$ from the right and it acts 
on a scheme $Y$ from the left and all schemes and actions are over 
$\Spec k$ for some field $k$, then 
we denote the quotient of $X\times_{\Spec k} Y$ by the diagonal $G$-action
defined by $(x,y)\mapsto (xg,g^{-1}y)$ for all $g\in G$ by
$$
   X\wedge^G Y \,:=\, (X\times_{\Spec k} Y)/G.
$$
We refer to \cite[Chapter III.1.3]{Giraud} for details and applications.

\begin{Proposition}
 \label{prop: product type classification}
 Let $X$ and $X\subset P$ be as in Theorem \ref{thm: product type}.
 \begin{enumerate}
      \item if $\rho(X)=2$, then
       $$\begin{array}{ccc}
         X &\iso& P'\times P'',
         \end{array}
       $$ 
       where $P'$ and $P''$ are Brauer--Severi curves over $k$, whose
       Brauer classes satisfy $[P]=[P']+[P'']\in\Br(k)$.
       In particular, $P\iso\PP^3_k$ if and only if $P'\iso P''$.
     \item If $\rho(X)=1$, then there exist a Brauer--Severi curve $P'$ over $k$
       and a finite Galois extension $K/k$ with Galois group $H:=\ZZ/2\ZZ$, 
       such that $X$ arises as twisted self-product
        $$\begin{array}{ccccc}
          X&\iso& (P'\times P')_K/H & = &\Spec K\wedge^H (P'\times P'),
          \end{array}
        $$
        where the $H$-action permutes the factors of $P'_K\times P'_K$.
        Moreover, $P\iso\PP^3_k$
        and $P'$ is a hyperplane section of $X\subset \PP^3_k$.
  \end{enumerate}
\end{Proposition}

\proof
We keep the notations and assumptions from the proof of 
Theorem \ref{thm: product type}.
The $G_k$-action fixes the class $(1,1)$.
Since the $G_k$-action preserves the intersection pairing on $\Pic_{X/k}(k^{\rm sep})$,
it follows that $G_k$ acts on $\ZZ(1,-1)$ either trivially, or by sign changes.
We have $\rho(X)=2$ in the first case, and $\rho(X)=1$ in the latter.

First, assume that $\rho(X)=2$.
By Theorem \ref{thm: main}, the classes
$(1,0)$ and $(0,1)$ give rise to morphisms
to Brauer--Severi curves $X\to P'$ and $X\to P''$ of class
$[P']=\delta((1,0))$ and $[P'']=\delta((0,1))$ in $\Br(k)$, respectively.
Thus, we obtain a morphism $X\to P'\times P''$, which is an isomorphism because
it is an isomorphism over $k^{\rm sep}$.
Since $\delta$ is a homomorphism, we find
$[P]=\delta(\calL)=\delta((1,1))=\delta((1,0))+\delta((0,1))=[P']+[P'']$.
Using that $P'$ and $P''$ are of period $2$, we find
that $P\iso\PP^3_k$ if and only if $[P]=0$, i.e., if and only if
$[P']=[P'']$.
By Corollary \ref{cor: isomorphic BS},
the latter is equivalent to $P'\iso P''$.

Second, assume that $\rho(X)=1$.
Then, the $G_k$-action permutes 
$(0,1)$ and $(1,0)$, i.e., it permutes the factors of
$\PP^1_{k^{\rm sep}}\times\PP^1_{k^{\rm sep}}$.
Thus, there exists a unique quadratic Galois extension $K/k$,
such that ${\rm Gal}(k^{\rm sep}/K)$ 
acts trivially on $\Pic_{X/k}(k^{\rm sep})$ and by the previous analysis
we have  $X_K:=Q''\times Q'''$ for two
Brauer--Severi curves $Q''$, $Q'''$ over $K$.
Using these and the $H:={\rm Gal}(K/k)$-action, we obtain a $H$-stable
diagonal embedding $Q'\subset X_K$ of a Brauer--Severi curve 
over $K$, and then, the two projections induce isomorphisms 
$Q'\iso Q''$ and $Q'\iso Q'''$ over $K$.
Taking the quotient by $H$, we obtain a Brauer--Severi curve
$P':=Q'/H\subset X$ over $k$.
Clearly, $P'_K\iso Q'$ and we obtain the description of
$X$ as twisted self-product.
On $X$, the curve $P'$ is a section of the class $(1,1)$,
which implies that this class comes from an invertible
sheaf, and thus, $0=\delta((1,1))\in\Br(k)$ by Proposition \ref{prop: delta}.
Since $\delta((1,1))=[P]$, we conclude $P\iso\PP^3_k$.
\qed\medskip

\begin{Remark}
  In the case of quadrics in $\PP^3$, similar results were already 
  established in \cite{CTSk}.
  A related, but somewhat different view on degree $8$ del~Pezzo surfaces 
  of product type was taken in (the proof of) \cite[Proposition 5.2]{CTKMdP6}:
  If $X$ is such a surface, then there exists a quadratic Galois 
  extension $K/k$ and a Brauer--Severi curve $C$ over $K$, 
  such that  $X\iso{\rm Res}_{K/k}C$,
  where ${\rm Res}_{K/k}$ denotes Weil restriction, see
  also \cite{Poonen}.
\end{Remark}

\begin{Corollary}
  \label{cor: h1 product type}
  Let $X$ be as in Theorem \ref{thm: product type}.
  Then,
  $$
     H^1\left(H,\, \Pic_{X/k}({k^{\rm sep}})\right) \,=\,0
  $$
  for all closed subgroups $H\subseteq G_k$, and
  $$
   {\rm Am}(X) \,\iso\, \left\{
   \begin{array}{ll}
     0 & \mbox{ if $\rho=1$ or if $X\iso\PP^1_k\times\PP^1_k$,}\\
     (\ZZ/2\ZZ)^2 & \mbox{ if $\rho=2$ and $\PP^1_k\not\iso P'\not\iso P''\not\iso\PP^1_k$,}\\
     (\ZZ/2\ZZ) & \mbox{ in the remaining $\rho=2$-cases.}
   \end{array}
   \right.
  $$
\end{Corollary}

\prf
Set $H^1(H):=H^1(H,\, \Pic_{X/k}({k^{\rm sep}}))$.
If $\rho=2$, then the $G_k$-action on $\Pic_{X/k}(k^{\rm sep})$ is trivial,
and we find $H^1(H)=0$ as in the proof of Lemma \ref{lem: h1 brauer-severi}.
Moreover, ${\rm Am}(X)$ is generated by $\delta((0,1)$ and $\delta((1,0))$, i.e.,
by $[P']$ and $[P'']$ in $\Br(k)$.
From this, the assertions on ${\rm Am}(X)$ follow in case $\rho=2$.

If $\rho=1$, then there exists an isomorphism $\Pic_{X/k}(k^{\rm sep})\iso\ZZ^2$, such that
the $G_k$-action factors through 
a surjective homomorphism $G_k\to\ZZ/2\ZZ$ and acts
on $\ZZ^2$ via $(a,b)\mapsto(b,a)$.
In particular, we find $H^1(\ZZ/2\ZZ,\ZZ^2)=0$ with respect to this action, see, 
for example, \cite[Chapter III.1, Example 2]{Brown}.
From this, we deduce $H^1(H)=0$ using inflation maps.
Moreover, ${\rm Am}(X)$ is generated by $\delta((1,1))$,
which is zero, since $(1,1)$ is the class of an invertible sheaf.
\qed\medskip

\begin{Corollary}
  \label{cor: product case point}
  If $X$ is as in Theorem \ref{thm: product type}, then
  the following are equivalent
  \begin{enumerate}
   \item $X$ is birationally equivalent to a Brauer--Severi surface,
   \item $X$ is rational,
   \item $X$ has a $k$-rational point, and
   \item $X$ is isomorphic to
    $$
       X\,\iso\,\PP^1_k\times \PP^1_k \mbox{ \quad or to \quad }X\,\iso\,\Spec K\wedge(\PP^1_k\times\PP^1_k).
    $$   
  \end{enumerate}
\end{Corollary}

\prf
The implications $(2)\Rightarrow(1)$ and $(2)\Rightarrow(3)$ are trivial,
and we established $(3)\Rightarrow(2)$ in Theorem \ref{thm: product type}.
Moreover, if $X$ is birationally equivalent to a Brauer--Severi surface $P$, 
then ${\rm Am}(P)={\rm Am}(X)$ is cyclic of order $1$ or $3$ 
by Lemma \ref{lem: h1 brauer-severi} and Theorem \ref{thm: birational invariance of h1}.
Together with Corollary \ref{cor: h1 product type}, we conclude
${\rm Am}(P)={\rm Am}(X)=0$, i.e.,
$P\iso\PP^2_k$, which establishes $(1)\Rightarrow(2)$.

Since $(4)\Rightarrow(3)$ is trivial, it remains to establish $(3)\Rightarrow(4)$.
Thus, we assume $X(k)\neq\emptyset$.
If $\rho=2$, then $X\iso P'\times P''$ and both Brauer--Severi curves $P'$ and $P''$
have $k$-rational points, i.e., $X\iso\PP^1_k\times\PP^1_k$.
If $\rho=1$, we have an embedding $X\subset\PP^3_k$
and $X\iso\Spec K\wedge (P'\times P')$.
Since $X(k)\neq\emptyset$, we have $X(K)\neq\emptyset$,
which yields $P'(K)\neq\emptyset$, and thus $P'_K\iso\PP^1_K$.
A $k$-rational point on $X$ gives rise to a $K$-rational and
${\rm Gal}(K/k)$-stable point on $X_K\iso\PP^1_K\times\PP^1_K$.
In particular, this point lies on some diagonal $\PP^1_K\subset X_K$,
and thus, lies on some diagonal $P''\subseteq X$ with
$X\iso\Spec K\wedge(P''\times P'')$.
Since $P''(k)\neq\emptyset$, we find $P''\iso\PP^1_k$.
\qed\medskip

We refer to Section \ref{sec: del Pezzo application} for more applications
of these results to the arithmetic and geometry of these surfaces.
 
\section{Del~Pezzo surfaces of large degree}
\label{sec: dP large degree}

Let $X$ be a del Pezzo surface of degree $d$ over a field $k$ that is not of product type.
Then, there exists a birational morphism
$$
 \begin{array}{ccccc}
   \overline{f} &:& \overline{X} &\to&\PP^2_{\overline{k}}
 \end{array}
$$
that is a blow-up in $(9-d)$ closed points $P_1,...,P_{9-d}$ in general position.
We set $H:=\overline{f}^*\OO_{\PP^2_{\overline{k}}}(1)$ and let $E_i:=\overline{f}^{-1}(P_i)$
be the exceptional divisors of $\overline{f}$.
Then, there exists an isomorphism of abelian groups
$$
 \begin{array}{ccc}
  \Pic(\overline{X}) &\cong& \ZZ H\,\oplus\,\bigoplus_{i=1}^{9-d}\,\ZZ E_i.
 \end{array}
$$
The $(-1)$-curves of $\overline{X}$ consist of the $E_i$, 
of preimages under $\overline{f}$ of lines through two distinct points $P_i$, 
of preimages under $\overline{f}$ of quadrics through five distinct points $P_i$, etc.,
and we refer to \cite[Theorem 26.2]{Manin} for details.
Let $K_{\overline{X}}$ be the canonical divisor class of $\overline{X}$, and let $\widetilde{E}$ 
be the sum of all $(-1)$-curves on $\overline{X}$.
We leave it to the reader to verify the following table. 
\begin{center}
  \begin{tabular}{|c||rcr|rcrcr|}\hline
    $d$ & \multicolumn{3}{c|}{class of $\widetilde{E}$ in $\Pic(\overline{X})$}  & 
    \multicolumn{5}{c|}{relations} 
    \\ \hline \hline
    $9$ & \multicolumn{3}{c|}{$0$}                 &  $3H$&$=$&$-K_{\overline{X}}$ &&\\
    $8$ & & & $E_1$                                      & $3H$&$=$&$-K_{\overline{X}}$&+&$\widetilde{E}$ \\
    $7$ & $H$ & &                                          & $H$&$=$& & & $\widetilde{E}$ \\
    \hline
    $6$ & $3H$&$-$&$\sum_{i=1}^3E_i$       & $0$&$=$&$-K_{\overline{X}}$ &-&$\widetilde{E}$ \\
    $5$ & $6H$&$-$&$2\sum_{i=1}^4E_i$     & $0$&$=$&$-2K_{\overline{X}}$&-&$\widetilde{E}$ \\
    $4$ & $12H$&$-$&$4\sum_{i=1}^5E_i$   & $0$&$=$&$-4K_{\overline{X}}$ &-&$\widetilde{E}$\\
    $3$ & $27H$&$-$&$9\sum_{i=1}^6E_i$   & $0$&$=$&$-9K_{\overline{X}}$ &-&$\widetilde{E}$\\
    $2$ & $84H$&$-$&$28\sum_{i=1}^7E_i$ & $0$&$=$&$-28K_{\overline{X}}$  &-&$\widetilde{E}$\\
    $1$ & $720H$&$-$&$240\sum_{i=1}^8E_i$   & $0$&$=$&$-240 K_{\overline{X}}$ &-&$\widetilde{E}$
   \\ \hline
\end{tabular}
\end{center}
Together with Theorem \ref{thm: main}, we obtain the following result.

 \begin{Theorem}
  \label{thm: del Pezzo descent}
 Let $X$ be a del~Pezzo surface of degree $d\geq7$ over a field $k$ that is
 not of product type.
 Then, $\overline{f}$ descends to a birational morphism
 $$
    f : X\to P
 $$
 to a Brauer--Severi surface $P$ over $k$, where
  $$
      \delta(H) \,=\, [P] \,\in\,\Br(k) \mbox{ \quad and \quad }{\rm Am}(X)\,\iso\,\ZZ/{\rm per}(P)\ZZ.
  $$
  Moreover, $X$ is rational if and only if $P\iso\PP^2_k$.
  This is equivalent to $X$ having a $k$-rational point.
 \end{Theorem}
 
\proof
By Theorem \ref{thm: Coombes}, the invertible sheaf $H$ on $X_{\overline{k}}$
defining $\overline{f}$ already lies in $\Pic_X({k^{\rm sep}})$, i.e.,
$\overline{f}$ descends to $k^{\rm sep}$,
and by abuse of notation, we re-define $\overline{X}$ to be $X_{k^{\rm sep}}$.
Clearly, the canonical divisor class $K_{\overline{X}}$ is 
$G_k$-invariant, and since
$G_k$ permutes the $(-1)$-curves of $\overline{X}$, also the class
of $\widetilde{E}$ is $G_k$-invariant.
In particular, $K_{\overline{X}}$ and $\widetilde{E}$ define classes in 
$\Pic_{X/k}(k^{\rm sep})^{G_k}\iso\Pic_{(X/k)({\rm fppf})}(k)$.
If $d\geq7$, then the above table shows that there exist 
positive multiples of $H$ that are integral linear combinations of $K_{\overline{X}}$ 
and $\widetilde{E}$.
Thus, $H\in\Pic_{X/k}(k^{\rm sep})$ descends to a class in $\Pic_{(X/k)({\rm fppf})}(k)$.
By Theorem \ref{thm: main}, $\overline{f}$ descends
to a birational morphism $f:X\to P$, where $P$ is a Brauer--Severi surface
of class $\delta(H)\in\Br(k)$.
The assertion on ${\rm Am}(X)$ follows from 
Proposition \ref{prop: amitsur is birational invariant} and 
Theorem \ref{thm: amitsur birational}.

If $X$ has a $k$-rational point, then so has $P$, and then $P\iso\PP^2_k$.
Since $f$ is a birational morphism, $P\iso\PP^2_k$ implies that $X$ is rational.
And if $X$ is rational, then it has a $k$-rational point by Lemma \ref{lem: Lang}.
\qed\medskip

As an immediate consequence, we obtain rationality and the existence of $k$-rational points
in some cases.
 
\begin{Corollary}
  \label{cor: has rational point}
  Let $X$ be as in Theorem \ref{thm: del Pezzo descent}.
  If $d\in\{7,8\}$, then $X$ has a $k$-rational point and
  $\overline{f}$ descends to a birational morphism
 $f:X\to\PP^2_k$.
\end{Corollary}

\proof
By Theorem \ref{thm: del Pezzo descent}, there exists a birational morphism
$X\to P$ that is a blow-up in a closed subscheme $Z\subset P$ of length $(9-d)$.
By Corollary \ref{cor: zero cycle}, we have $P\iso\PP^2_k$ if $3$ and $(9-d)$ are coprime.
In particular, we have $X(k)\neq\emptyset$ in these cases by
Theorem \ref{thm: del Pezzo descent} and Lemma \ref{lem: Lang}.
\qed\medskip

Since a del~Pezzo surface of degree $9$ is a Brauer--Severi 
surface, it has rational points if and only if it is trivial.
In particular, Corollary \ref{cor: has rational point} does not hold
for $d=9$.
 
\subsection{Applications to arithmetic geometry}
\label{sec: del Pezzo application}
We now give a couple of applications of the just established results.
Again, we stress that most if not all of these applications are well-known, and 
merely illustrate the usefulness of studying varieties via Brauer--Severi
varieties.

\begin{Corollary}
  \label{cor: h1 for non product} 
  If $X$ is a del~Pezzo surface of degree $\geq7$ over $k$, then
  $$\begin{array}{ccc}
     H^1\left(H,\, \Pic_{X/k}({k^{\rm sep}})\right) &=&0.
    \end{array}
  $$
  for all closed subgroups $H\subseteq G_k$
 \end{Corollary}
 
\prf
If $X$ is not of product type, then it is birationally equivalent to a Brauer--Severi surface $P$
by Theorem \ref{thm: del Pezzo descent}, and then
the statement follows from Theorem \ref{thm: birational invariance of h1}
and Lemma \ref{lem: h1 brauer-severi}. 
If $X$ is of product type, then this is Corollary \ref{cor: h1 product type}.
\qed\medskip

For the next application, let us recall that a surface is called {\em rational} 
if it is birationally equivalent to $\PP^2_k$, and that it is called {\em unirational}
if there exists a dominant and rational map from $\PP^2_k$ onto it.
The following result is a special case of \cite[Theorem 29.4]{Manin}.
 
\begin{Corollary}
 \label{cor: birational to p2}
 Let $X$ be a del~Pezzo surface of degree $\geq7$ over a field $k$.
 Then, the following are equivalent:
 \begin{enumerate}
  \item $X$ is rational,
  \item $X$ is unirational, and
  \item $X$ has a $k$-rational point.
 \end{enumerate}
\end{Corollary}

\prf
Clearly, we have $(1)\Rightarrow(2)\Rightarrow(3)$, whereas $(3)\Rightarrow(1)$ 
follows from Corollary \ref{cor: product case point} and 
Theorem \ref{thm: del Pezzo descent}.
\qed\medskip

This leads us to the question whether a del~Pezzo surface necessarily has 
a $k$-rational point.
Over finite fields, this is true and follows from the Weil conjectures, which we will
recall in Theorem \ref{thm: weil} below.
By a theorem of Wedderburn, finite fields have trivial Brauer groups, 
and thus, the following corollary gives existence of $k$-rational points for
more general fields.

 \begin{Corollary}
   Let $X$ be a del~Pezzo surface of degree $\geq7$ over a field $k$ with 
   $\Br(k)=0$.
   Then, $X$ has a $k$-rational point, and thus, is rational.
 \end{Corollary}
 
 \proof
If $X$ is not of product type, then there exists a birational morphism
$f:X\to P$ to a Brauer--Severi surface by Theorem \ref{thm: del Pezzo descent}.
Since $\Br(k)=0$, we have $P\iso\PP^2_k$, and Theorem \ref{thm: del Pezzo descent}
gives $X(k)\neq\emptyset$.
 
Thus, let $X$ be of product type.
By Proposition \ref{prop: product type classification},
$X$ is a product of Brauer--Severi curves ($\rho=2$), 
or contains at least a Brauer--Severi curve ($\rho=1$).
Since $\Br(k)=0$, all Brauer--Severi curves are isomorphic to $\PP^1_k$,
and thus, contain $k$-rational points.
In particular, we find $X(k)\neq\emptyset$.
\qed\medskip

In Section \ref{subsec: Hasse}, we discussed the Hasse principle
and weak approximation for varieties over global fields.
Here, we establish the following.

 \begin{Corollary}
   \label{cor:hasse principle}
   Del~Pezzo surfaces of degree $\geq7$ over global fields satisfy 
   weak approximation and the Hasse principle.
 \end{Corollary}
 
\proof
If $X$ is not of product type, then it is birationally equivalent to a Brauer--Severi surface by
Theorem \ref{thm: del Pezzo descent}, and since the two claimed properties are
preserved under birational maps and hold for Brauer--Severi varieties,
the assertion follows in this case.

If $X$ is of product type, then there are two 
cases by Proposition \ref{prop: product type classification}.
If $\rho=2$, then $X$ is a product of two Brauer--Severi curves, and we conclude as before.

Thus, we may assume $\rho=1$.
Let us first establish the Hasse principle:
there exists a quadratic Galois extension $L/K$, such that $\rho(X_L)=2$.
From $X(K_\nu)\neq\emptyset$ for all $\nu\in\Omega_K$, we find
$X_{L_\mu}\iso\PP^1_{L_\mu}\times\PP^1_{L_\mu}$ for all $\mu\in\Omega_L$, 
and thus, $X_L\iso\PP^1_L\times\PP^1_L$ by the Hasse principle for
Brauer--Severi curves.
As in the proof of Corollary \ref{cor: product case point}, 
we exhibit $X$ as twisted self-product of $\PP^1_k$, which has a $k$-rational
point and establishes the Hasse principle.
Thus, to establish weak approximation, we may assume that $X$ has
a $k$-rational point.
But then, $X$ is rational by Corollary \ref{cor: product case point},
and since weak approximation is a birational invariant, the assertion follows.
\qed\medskip

\section{Del~Pezzo surfaces of degree 6}

In the previous sections, we have seen a close connection between Brauer--Severi varieties
and del~Pezzo surfaces of  degree $\geq7$.
In this section, we discuss del~Pezzo surfaces of degree $6$, which are not so 
directly linked to Brauer--Severi varieties.

For the geometry and the arithmetic of these surfaces, we refer the 
interested reader to \cite{CTdP6}, \cite{Manin}, and the survey \cite[Section 2.4]{Varilly}.
We keep the notation introduced in Section \ref{sec: dP large degree}:
If $X$ is a degree $6$ del~Pezzo surface over a field $k$, then
there exists a blow-up $f_{\overline{k}}:\overline{X}\to\PP^2_{\overline{k}}$ in
three points in general position
with exceptional $(-1)$-curves $E_1$, $E_2$, and $E_3$.
Then, there are six $(-1)$-curves on $X$, namely the three exceptional curves
$E_i$, $i=1,2,3$ of $\overline{f}$, as well as 
the three curves $E_i':=H-E_j-E_k$, $i=1,2,3$ where $\{j,k\}=\{1,2,3\}\backslash\{i\}$
and where $H=\overline{f}^*\OO_{\PP^2}(1)$ as in Section \ref{sec: dP large degree}.
These curves intersect in a hexagon as follows.
$$
 
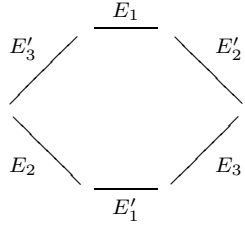
\xymatrix{ 
   &  \ar@{-}^{E_1}[r] & \ar@{-}^{E_2'}[dr]  &\\
    \ar@{-}^{E_3'}[ur]\ar@{-}_{E_2}[dr] & & &  \\
    &  \ar@{-}_{E_1'}[r] & \ar@{-}_{E_3}[ur]  &
  }
$$
The absolute Galois group $G_k$ acts on these six $(-1)$-curves on $X_{k^{\rm sep}}$,
and associated to this action, we have following field extensions of $k$.
\begin{enumerate}
  \item Since $G_k$ acts on the two sets $\{E_1,E_2,E_3\}$ and $\{E_1',E_2',E_3'\}$,
    there is a group homomorphism
    $$
       \varphi_1\,:\,G_k\,\to\,S_2\,\iso\,\ZZ/2\ZZ.
    $$
    The fixed field of either of the two sets is a finite separable extension $k\subseteq K$ with $[K:k]|2$,
    and $k\neq K$ if and only if $\varphi_1$ is surjective.
  \item Since $G_k$ acts on the three sets $\{E_i,E_i'\}$, $i=1,2,3$, 
    there is a group homomorphism
    $$
      \varphi_2\,:\,G_k\,\to\,S_3.
    $$
    There exists a finite separable extension $k\subseteq L$ with $[L:k]|3$, unique up to conjugation
    in $k^{\rm sep}$, over which at least one of these three sets is defined.
    We have $k\neq L$ if and only if $3$ divides the order of $\varphi_2(G_k)$.
    Next, there exists a finite and separable extension
    $L\subseteq M$ with $[M:L]|2$, over which all three sets are defined.
\end{enumerate}

Combining $\varphi_1$ and $\varphi_2$, we obtain a group homomorphism
$$
   G_k\,\stackrel{\varphi_1\times\varphi_2}{\longrightarrow}\, \ZZ/2\ZZ\times S_3\,\iso\,D_{2\cdot 6},
$$
where $D_{2\cdot 6}$ denotes the dihedral group of order $12$, i.e., the automorphism group
of the hexagon.
Using these field extensions, we obtain the following classification, which uses and 
slightly extends a classical result of Manin from \cite{Manin} in case (3).

\begin{Theorem}
  \label{thm: degree 6}
  Let $X$ be a del~Pezzo surface of degree $6$ over a field $k$.
  \begin{enumerate}
   \item The morphism $\overline{f}$ descends to a birational morphism 
     $$ 
        f\,:\,X\,\to\,P
     $$
     to a Brauer--Severi surface $P$ if and only if $k=K$.
     In this case,  $\rho(X)\geq2$ and ${\rm Am}(X)={\rm Am}(P)$.
   \item There exists a birational morphism $X\to Y$
     onto a degree $8$ del~Pezzo surface $Y$ of product type if and only if
     $k=L$.
     In this case, 
     $$\begin{array}{c|ccc}
        & \rho(X) & Y  \\
       \hline
       k\neq M & 3 & \Spec M\wedge(\PP^1_k\times\PP^1_k) \\       
       k=M         & 4 & \PP^1_k\times\PP^1_k  
     \end{array}$$
     $X$ has a $k$-rational point, and ${\rm Am}(X)=0$.
   \item If $k\neq K$ and $k\neq L$, then $\rho(X)=1$, ${\rm Am}(X)=0$, and the following are equivalent.
    \begin{enumerate}
     \item $X$ is birationally equivalent to a Brauer--Severi surface,
     \item $X$ is birationally equivalent to a product of two Brauer--Severi curves,
     \item $X$ is rational, and
     \item $X$ has a $k$-rational point.
    \end{enumerate}
  \end{enumerate}
\end{Theorem}

\prf
Let us first show (1).
If $k=K$, then $F:=E_1+E_2+E_3$ descends to a class in $\Pic(X_{k^{\rm sep}})^{G_k}=\Pic_{(X/k)({\rm fppf})}(k)$
and  we find $\rho(X)\geq2$.
Thus, also $H=\frac{1}{3}(-K_X+F)$ descends to a class in $\Pic_{(X/k)({\rm fppf})}(k)$, and by
Theorem \ref{thm: main}, we obtain a birational morphism $|H|:X\to P$ to a Brauer--Severi surface,
which coincides with $\overline{f}$ over $\overline{k}$.
Conversely, if $\overline{f}$ descends to a birational morphism $f:X\to P$,
then the exceptional divisor of $f$ is of class $F$ or $E'_1+E_2'+E_3'$, 
and we find $k=K$.
Moreover, we have ${\rm Am}(X)={\rm Am}(P)$ by 
Theorem \ref{thm: birational invariance of h1}.

If $k=L$, then, say $E_1+E_1'$, descends to a class in $\Pic(X_{k^{\rm sep}})^{G_k}$.
Moreover, we find that the classes
$\frac{1}{2}(-K_X+E_1+E_1')=\,2H-E_2-E_3$ as well as 
$\frac{1}{2}(-K_X-E_1-E_1')\,=\, H-E_1$,
and thus, the classes $H$, $E_1$, and $E_1'=H-E_2-E_3$ lie in 
$\Pic(X_{k^{\rm sep}})^{G_k}$.
The $G_k$-action is trivial on $H$ and $E_1$, whereas
it is either trivial on the set $\{E_2,E_3\}$ (if $k=M$)
or permutes the two (if $k\neq M$).
Since the class of $E_1$ is $G_k$-invariant and there is a unique effective
divisor in this linear system, we find that $\PP^1_k\iso E_1\subset X$.
In particular, $X$ has a $k$-rational point and ${\rm Am}(X)=0$.
Using Theorem \ref{thm: main} and the fact that $X$ has a $k$-rational point,
we obtain a birational morphism
$$
   |\frac{1}{2}(-K_X+E_1+E_1')|\,:\, X\to Y\subset \PP^3_k
$$
onto a smooth quadric $Y$ with a $k$-rational point.
In particular, $Y$ is a degree $8$ del~Pezzo surface of product type.
Over $k^{\rm sep}$, this morphism contracts $E_1$ and $E_1'$ and thus,
we find
$$
  \Pic(Y_{k^{\rm sep}}) \,\iso\, \left( \ZZ H\oplus\bigoplus_{i=1}^3\ZZ E_i\right)/
  \langle E_1, E_1' \rangle
  \,\iso\,\ZZ \overline{E}_2 \oplus \ZZ \overline{E}_3.
$$
The $G_k$-action on it is either trivial ($k=M$) or permutes the two summands ($k\neq M$).
Using $Y(k)\neq\emptyset$ and Corollary \ref{cor: product case point},
we find $\rho(X)=4$ and $Y\iso\PP^1_k\times\PP^1_k$ in the first case,
and $\rho(X)=3$ and $Y\iso\Spec M\wedge(\PP^1_k\times\PP^1_k)$ in 
the latter.
Conversely, if there exists a birational morphism $X\to Y$ onto a degree $8$ del~Pezzo
surface of product type, then the exceptional divisor is of class 
$E_i+E_i'$ for some $i$, and thus, $k=L$.
This establishes (2).

Finally, assume that $k\neq K$ and $k\neq L$.
Then, $\varphi_1$ is surjective, and $\varphi_2(G_k)$ contains all $3$-cycles of $S_3$.
From this, it is not difficult to see 
that $\Pic(X_{\overline{k}})^{G_k}$ is of rank $1$ and generated by the class of $K_X$.
Since this latter class is an invertible sheaf, we find ${\rm Am}(X)=0$.
Thus, if $X$ is birationally equivalent to a Brauer--Severi surface $P$,
then ${\rm Am}(X)=0$ together with Lemma \ref{lem: h1 brauer-severi}   
and Theorem \ref{thm: birational invariance of h1} implies that
$P\iso\PP^2_k$.
Similarly, if $X$ is birationally equivalent to the product $P'\times P''$ of two Brauer--Severi curves,
then $P'\iso P''\iso\PP^1_k$.
From this, we obtain the implications
$(a)\Leftrightarrow (b)\Leftrightarrow (c)\Rightarrow (d)$.
The implication $(d)\Rightarrow (c)$ is due to Manin \cite[Theorem 29.4]{Manin}.
\qed\medskip

\begin{Remark}
 \label{rem: amitsur remark}
  In case (1) of the above Theorem it is important to note that $P$ need not be unique, 
  but that ${\rm Am}(P)$ is well-defined.
  More precisely, if we set $F:=E_1+E_2+E_3$ and 
  $F'=E_1'+E_2'+E_3'$, then Theorem \ref{thm: main} provides us with two morphisms to 
  Brauer--Severi surfaces $P_1$ and $P_2$
  $$
  \begin{array}{lclcc}
     |H| &=& |\frac{1}{3}(-K_X+F)| &:& X \,\to\, P_1 \\
     |H'| &:=&|\frac{1}{3}(-K_X+F')|  &:& X\,\to\, P_2\\
  \end{array}
  $$
  Since $H+H'=-K_X$ and $\delta(K_X)=0$, we find 
  $$
    [P_1]\,=\,\delta(H)\,=\,\delta(-K_X-H')\,=\,-\delta(H')\,=\,-[P_2]\,\in\,\Br(k),
  $$
  and thus, $P_1\iso P_2$ if and only if both are isomorphic to $\PP^2_k$.
  On the other hand, $P_1$ and $P_2$ are birationally equivalent, since we have birational morphisms
  $$
     P_1\,\stackrel{|H|}{\longleftarrow}\,X\,\stackrel{|H'|}{\longrightarrow} P_2\,.
  $$
  Over $\overline{k}$, this becomes the blow-up of three closed points $Z$ followed by the blow-down
  of the three $(-1)$-curves that are the strict transforms of lines through any two of the points in $Z$.
   This is an example of a {\em Cremona transformation}.  
\end{Remark}

We remark that a surface of case (3) and without $k$-rational points is
neither birationally equivalent to a Brauer--Severi surface nor to the product of two Brauer--Severi curves.
For finer and more detailed classification results for degree $6$ del~Pezzo surfaces, we refer the interested 
reader to \cite{Corn}, \cite{Blunk}, and \cite{CTKMdP6}.
Finally, the sum $\widetilde{E}$ of all $(-1)$-curves on $X_{k^{\rm sep}}$ is a $G_k$-invariant divisor, 
and thus, descends to a curve on $X$.
By \cite[Theorem 30.3.1]{Manin}, the complement $X\backslash\widetilde{E}$ is isomorphic to
a torsor under a two-dimensional torus over $k$, which can be used to study the arithmetic and geometry of 
these surfaces, see also \cite{Skorobogatov book}.

\section{Del~Pezzo surfaces of small degree}
\label{sec: small degree}

For the remainder of this article, our results will be less complete and less self-contained.
We will circle around questions of birationality of a del~Pezzo surface $X$
of degree $\leq5$ to Brauer--Severi surfaces, and about descending the morphism 
$\overline{f}:\overline{X}\to\PP^2_{\overline{k}}$ to $k$.

\subsection{Birationality to Brauer--Severi surfaces}
Let $k=\FF_q$ be a finite field of characteristic $p$, and let
$X$ be a smooth and projective surface over $k$ 
such that $X_{\overline{k}}$ is birationally equivalent to $\PP^2$.
Then, it follows from the Weil conjectures (in this case
already a theorem of Weil himself) that
the number of $k$-rational points is congruent to
$1$ modulo $q$, see \cite[Chapter IV.27]{Manin}.
In particular, we obtain that

\begin{Theorem}[Weil]
  \label{thm: weil}
  If $X$ is a del~Pezzo surface over a finite field $\FF_q$, then $X$
  has a $\FF_q$-rational point.
\end{Theorem}

Since ${\rm Br}(\FF_q)=0$ by a theorem of Wedderburn,
there are no non-trivial Brauer--Severi varieties over $\FF_q$.

\begin{Remark}
  \label{rem: M-CT-S}
  Let $X$ be a del~Pezzo surface of degree $\geq5$ over a field $k$.
  Manin \cite[Theorem 29.4]{Manin} showed that $X$ is
  rational if and only if it contains a $k$-rational point.
  Even if $X$ has no $k$-rational point,
  Manin \cite[Theorem 29.3]{Manin}
  showed that
  $$
    H^1\left(H,\,\Pic_{(X/k)({\rm fppf})}(k^{\rm sep})\right)\,=\,0
  $$ 
  for all closed subgroups $H\subseteq G_k$.
  We refer to \cite[Th\'eor\`eme 2.B.1]{CTS} for a general
  principle explaining this vanishing of cohomology.
\end{Remark}

In this section, we give a partial generalization to birational maps to
Brauer--Severi surfaces.

\begin{Lemma}
 \label{lem: del Pezzo Amitsur}
 Let $X$ be a degree $d$ del~Pezzo surface over $k$.
 Then, 
 \begin{enumerate}
  \item There exists an effective zero-cycle $Z$ of degree $d$ on $X$.
    If $d\neq2$ or if ${\rm char}(k)\neq2$, then there exists such a zero-cycle 
    $Z$, whose closed points have residue fields that are separable over $k$.
  \item The abelian group ${\rm Am}(X)$ is finite and every element has an order
    dividing $d$.
 \end{enumerate}
\end{Lemma}

\prf
If $d\geq3$, then $\omega_X^{-1}$ is very ample, and $|\omega_X^{-1}|$
embeds $X$ as a surface of degree $d$ into $\PP^d_k$.
Intersecting $X$ with a linear subspace of codimension $2$, 
we obtain an effective zero-cycle $Z$ of degree $d$ on $X$.
The closed points of $Z$ have automatically separable residue fields if
$k$ is finite.
Otherwise, $k$ is infinite, and then, the intersection with a generic linear
subspace of codimension $2$ yields a $Z$ that is smooth over $k$ by
 \cite[Th\'eor\`eme I.6.3]{Bertini theorems}.
Thus, in any case, we obtain a $Z$, whose closed points have residue fields
that are separable over $k$.
If $d=2$, then $|\omega_X^{-1}|$  defines a double cover $X\to\PP^2_k$,
and the pre-image of a $k$-rational point yields an effective
zero-cycle $Z$ of degree $2$ on $X$.
If ${\rm char}(k)\neq2$, then residue fields of closed points of $Z$
are separable over $k$.
If $d=1$, then $|-K_X|$ has a unique-base point, and in particular, 
$X$ has a $k$-rational point.
This establishes (1).
Since $b_1(X)=0$, the group ${\rm Am}(X)$ is finite by Lemma \ref{lem: picard rank}.
Then, assertion (2) follows from Corollary \ref{lem: order amitsur}.
\qed\medskip

\begin{Corollary}
  \label{cor: dP has a rational point}
   Let $X$ be a del~Pezzo surface of degree $d$ over a field $k$.
   \begin{enumerate}
     \item If $d\in\{1,2,4,5,7,8\}$ and $X$ is birationally equivalent to a Brauer--Severi surface $P$, then
       $P\iso\PP^2_k$ and $X$ has a $k$-rational point.
     \item If $d\in\{1,3,5,7,9\}$ and $X$ is birationally equivalent to a product $P'\times P''$ of two Brauer--Severi
        curves, then $P'\iso P''\iso\PP^1_k$ and $X$ has a $k$-rational point.
   \end{enumerate}
\end{Corollary}

\prf
Let $X$ and $d$ be as in (1).
Then, every element of ${\rm Am}(X)$ is of order dividing $d$ by Lemma \ref{lem: del Pezzo Amitsur},
but also of order dividing $3$ by Theorem \ref{thm: birational invariance of h1}
and Theorem \ref{thm: brauer severi picard}.
By our assumptions on $d$, we find ${\rm Am}(P)=0$, and thus, $P\iso\PP^2_k$.
Since the latter has a $k$-rational point, so has $X$ by Lemma \ref{lem: Lang}.
This shows (1).
The proof of (2) is similar and we leave it to the reader.
\qed\medskip

Combining this with a result of  Coray \cite{Coray}, we obtain the following.

\begin{Theorem}
  \label{thm: birational to BS surface}
  Let $X$ be a del~Pezzo surface of degree $d\in\{5,7,8\}$ over a perfect field $k$.
  Then, the following are equivalent
  \begin{enumerate}
   \item There exists a dominant and rational map $P\dashrightarrow X$ 
     from a Brauer--Severi surface $P$ over $k$,
   \item $X$ is birationally equivalent to a Brauer--Severi surface,
   \item $X$ is rational, and
   \item $X$ has a $k$-rational point.
  \end{enumerate}
\end{Theorem}

\prf
The implications $(3)\Rightarrow(2)\Rightarrow(1)$ are trivial.

Let $\varphi:P\dashrightarrow X$ be as in (1).
By Lemma \ref{lem: del Pezzo Amitsur}, there exists a zero-cycle of
degree $9$ on $P$, and another one of degree $d$ on $X$.
Using $\varphi$, we obtain a zero-cycle of degree dividing $9$
on $X$.
By assumption, $d$ is coprime to $9$, and thus, 
there exists a zero-cycle of degree $1$ on $X$.
By \cite{Coray}, this implies that $X$ has a $k$-rational point and 
establishes $(1)\Rightarrow(4)$.

The implication $(4)\Rightarrow(3)$ is a result of Manin \cite[Theorem 29.4]{Manin} .
\qed\medskip

Now, if a del~Pezzo surface $X$ over a field $k$ is birationally equivalent to a Brauer--Severi surface,
then $H^1(H,\Pic_{X/k}(k^{\rm sep}))=0$ for all closed subgroups $H\subseteq G_k$
by Theorem \ref{thm: birational invariance of h1}.
Moreover, this vanishing holds for all del~Pezzo surfaces of degree $\geq5$, see
Remark \ref{rem: M-CT-S}.
However, for del~Pezzo surfaces of degree $\leq4$, these cohomology groups may be non-zero,
see \cite[Section 31]{Manin}, \cite{Sir Peter}, \cite{Kresch Tschinkel}, and \cite{Varilly weak}.
In particular, del~Pezzo surfaces of degree $\leq4$ are in general {\em not} birationally equivalent
to Brauer--Severi surfaces.

For further information concerning geometrically rational surfaces, 
unirationality, central simple algebras, and connections with cohomological dimension, 
we refer the interested reader to \cite{CTKMdP6}.

\subsection{Del~Pezzo surfaces of degree 5}

In order to decide whether a birational map 
$f_{\overline{k}}:X_{\overline{k}}\to\PP^2_{\overline{k}}$ as in 
Section \ref{sec: dP large degree} descends to $k$ for a degree $5$ del~Pezzo surface $X$ over $k$, 
we introduce the following notion.

\begin{Definition}
 \label{def: conic}
  Let $X$ be a del~Pezzo surface over a field $k$.
  A {\em conic} on $X$ is a geometrically integral curve $C$ on $X$
  with $C^2=0$ and $-K_X\cdot C=2$.
  An element $\calL\in\Pic_{(X/k)({\rm fppf})}(k)$ is called a
  {\em conic class} if $\calL\otimes_k\overline{k}\iso\OO_{X_{\overline{k}}}(\overline{C})$ 
  for some conic $\overline{C}$ on $X_{\overline{k}}$.
\end{Definition}

The following is an analogue of Theorem \ref{thm: del Pezzo descent}
for degree $5$ del~Pezzo surfaces.

\begin{Theorem}
  \label{thm: dP5 descent}
  Let $X$ be a del~Pezzo surface of degree $5$ over a field $k$.
  Then, the following are equivalent:
  \begin{enumerate}
   \item There exists a birational morphism $f:X\to P$ to a Brauer--Severi surface,
     such that $f_{\overline{k}}$ is the blow-up of $4$ points in general position.
   \item There exists a birational morphism $f:X\to\PP^2_k$,
     such that $f_{\overline{k}}$ is the blow-up of $4$ points in general position.
   \item There exists a class $F\in\Pic_{(X/k)({\rm fppf})}(k)$ such that 
     $$
        F_{\overline{k}} \,\iso\, \OO_{\overline{X}}(E_1+E_2+E_3+E_4),
     $$
     where the $E_i$ are disjoint $(-1)$-curves on $\overline{X}$.
    \item  There exists a conic class in $\Pic_{(X/k)({\rm fppf})}(k)$.
 \end{enumerate}
 If these equivalent conditions hold, then $X$ has a $k$-rational point.
\end{Theorem}

\prf
If $f$ is as in (1), then $X$ has a $k$-rational point by Corollary \ref{cor: dP has a rational point}.
Thus, $P\iso\PP^2_k$, and we obtain $(1)\Rightarrow(2)$.

If $f$ is as in $(2)$, then the exceptional divisor of $f$ is a class $F$ as stated in $(3)$,
and we obtain $(2)\Rightarrow(3)$.

If $f$ is as in $(3)$, then, using Theorem \ref{thm: main}, there exists a birational morphism
$|\frac{1}{3}(-K_X-F)|:X\to P$ to a Brauer--Severi surface $P$ as in $(1)$, which establishes
$(3)\Rightarrow(1)$.

If $f$ is as in $(2)$, let $Z\subset\PP^2_k$ be the degree $4$ cycle blown up by $f$.
Then $f^*(\OO_{\PP^2_k}(2)(-Z))$, i.e., the pullback of the pencil of conics through $Z$,
is a conic class on $X$ and establishes $(2)\Rightarrow(4)$.

Finally, if $C$ is a conic class on $X$, then, using Theorem \ref{thm: main}, 
there exists a birational morphism
$|-K_X+C|:X\to P$ to a Brauer--Severi surface $P$ as in $(1)$, which establishes
$(4)\Rightarrow(1)$.
\qed\medskip

\begin{Remark}
  By theorems of Enriques, Swinnerton-Dyer, Skorobogatov, Shepherd-Barron,
  Koll\'ar, and Hassett (see \cite[Theorem 2.5]{Varilly} for precise references and overview),
  a degree $5$ del~Pezzo $X$ over a field $k$ always has a $k$-rational point.
  Thus, $X$ is rational by \cite[Theorem 29.4]{Manin}, and we have
  $$
     {\rm Am}(X)=0,\mbox{ \quad as well as \quad } H^1(H,\Pic_{X/k}(k^{\rm sep}))=0
  $$
  for every closed subgroup $H\subseteq G_k$ by Corollary \ref{cor: amitsur trivial},
  Theorem \ref{thm: birational invariance of h1}, and Lemma \ref{lem: h1 brauer-severi}.
\end{Remark}

\subsection{Del~Pezzo surfaces of degree 4}
A classical theorem of Manin \cite[Theorem 29.4]{Manin} states
that a del~Pezzo surface of degree $4$ over a sufficiently large 
field $k$ is unirational if and only if it contains a $k$-rational point.
Here, we have the following analogue in our setting.

\begin{Proposition}
  Let $X$ be a del~Pezzo surface of degree $4$ over a perfect field $k$.
  Then, the following are equivalent
  \begin{enumerate}
    \item There exists a dominant rational map $P\dashrightarrow X$ from
      a Brauer--Severi surface $P$ over $k$.
   \item $X$ is unirational,
   \item $X$ has a $k$-rational point,
  \end{enumerate} 
\end{Proposition}

\prf
The implications $(2)\Rightarrow(1)$ is trivial and
$(2)\Rightarrow(3)$ is Lemma \ref{lem: Lang}.

The implication $(3)\Rightarrow(2)$ is shown in \cite[Theorem 29.4]{Manin}
and \cite[Theorem 30.1]{Manin} if $k$ has at least 23 elements
and in \cite[Theorem 2.1]{Knecht} and 
\cite[Proposition 5.19]{Pieropan} in the remaining cases.

To show $(1)\Rightarrow(3)$, we argue as in the proof of
the implication $(1)\Rightarrow(4)$ of Theorem \ref{thm: birational to BS surface}
by first exhibiting a degree $1$ zero-cycle on $X$, and then, 
using \cite{Coray} to deduce the existence of a
$k$-rational point on $X$.
We leave the details to the reader.
\qed\medskip

If a field $k$ is finite or perfect of characteristic $2$, then 
a degree $4$ del~Pezzo surface over $k$ always has a $k$-rational
point, see \cite[Theorem 27.1]{Manin} and \cite{Duncan}.
In this case, we also have ${\rm Am}(X)=0$.
From Lemma \ref{lem: del Pezzo Amitsur}, we infer that ${\rm Am}(X)$ 
is at most $4$-torsion for degree $4$ del~Pezzo surfaces.
For the possibilities of $H^1(G_k,\Pic_{X/k}(k^{\rm sep}))$, see \cite{Sir Peter}.

The following is an anolog of Theorem \ref{thm: del Pezzo descent}
for degree $4$ del~Pezzo surfaces.

\begin{Theorem}
  \label{thm: dP4 quasi-split}
  Let $X$ be a del~Pezzo surface of degree $4$ over a field $k$.
  Then, the following are equivalent:
  \begin{enumerate}
   \item There exists a birational morphism $f:X\to P$ to a Brauer--Severi surface,
     such that $f_{\overline{k}}$ is the blow-up of $5$ points in general position.
   \item There exists a birational morphism $f:X\to\PP^2_k$,
     such that $f_{\overline{k}}$ is the blow-up of $5$ points in general position.
   \item There exists a curve $\PP^1_k\iso E\subset X$ with $E^2=-1$.
   \item There exists a class $E\in\Pic_{(X/k)({\rm fppf})}(k)$ with $E^2=K_X\cdot E=-1$.
 \end{enumerate}
 If these equivalent conditions hold, then $X$ has a $k$-rational point.
\end{Theorem}

\prf
The implication $(2)\Rightarrow(1)$ is trivial.
If $f$ is as in (1), then $X$ has a $k$-rational point by Corollary \ref{cor: dP has a rational point}.
Thus, $P\iso\PP^2_k$, and we obtain $(1)\Rightarrow(2)$.

If $f$ is as in $(2)$, let $Z\subset\PP^2_k$ be the degree $5$ cycle blown up by $f$.
Then $f^*(\OO_{\PP^2_k}(2)(-Z))$, i.e., the pullback of the class of the unique conic 
through $Z$, is a class $E$ as stated in $(4)$ on $X$ and establishes 
$(2)\Rightarrow(4)$.

If $E$ is a class as in $(4)$, then, using Theorem \ref{thm: main}, there exists a birational
morphism $|-K_X-E|:X\to P$ to a Brauer--Severi surface $P$ as in $(1)$, which establishes
$(4)\Rightarrow(1)$.

The implication $(3)\Rightarrow(4)$ is trivial, and if $E$ is a class as in $(4)$,
then there exists a unique section of the associated invertible sheaf on $k^{\rm sep}$.
This is necessarily $G_k$-invariant, thus, descends to a curve on $X$,
and establishes $(4)\Rightarrow(3)$.
\qed\medskip

\begin{Remark}
In \cite{Skorobogatov}, Skorobogatov called del~Pezzo surfaces of degree $4$
that satisfy condition (3) above {\em quasi-split}.
\end{Remark}

Before proceeding, let us recall a couple of classical results on the geometry
of degree $4$ del~Pezzo surfaces,
and refer the interested reader to \cite{Skorobogatov} and 
\cite[Chapter 8.6]{Dolgachev} for details.
The anti-canonical linear system embeds $X$ 
as a complete intersection of two quadrics in $\PP^4_k$,
i.e., $X$ is given by $Q_0=Q_1=0$, where $Q_0$ and $Q_1$
are two quadratic forms in five variables over $k$.
The {\em degeneracy locus} of this pencil of quadrics
$$
 {\rm Deg}_X \,:=\, \left\{\,\det(t_0Q_0+t_1Q_1)\,=\,0 \,\right\} \,\subset\,
 \PP^1_k\,=\,\Proj k[t_0,t_1]
$$
is a zero-dimensional subscheme, which is \'etale and of length $5$ over $k$.
Over $\overline{k}$, its points correspond to the singular quadrics containing $X$, 
all of which are cones over smooth quadric surfaces. 
Let $\nu_2:\PP^1_k\to\PP^2_k$ be the $2$-uple Veronese embedding and set
$$
Z := \nu_2({\rm Deg}_X)\,\subset\,C:=\nu_2(\PP^1_k)\,\subset\,\PP^2_k.
$$
If $X$ contains a $k$-rational $(-1)$-curve, i.e., 
if $X$ is quasi-split, then $X$ is the blow-up of $\PP^2_k$ in $Z$, 
see Theorem \ref{thm: dP4 quasi-split} and \cite[Theorem 2.3]{Skorobogatov}.

\begin{Proposition}
 \label{prop: dP4 to dP8}
 Let $X$ be a del~Pezzo surface of degree $4$ over a field 
 $k$ of characteristic $\neq2$ with at least $5$ elements.
 Then, the following are equivalent:
 \begin{enumerate}
   \item The degeneracy scheme ${\rm Deg}_X$ has a $k$-rational point.
   \item There exists a finite morphism $\psi:X\to S$ of degree $2$, 
     where $S$ is a del~Pezzo surface of degree $8$ of product type.
 \end{enumerate}
 Moreover, if $\psi$ is as in $(2)$, then $S$ is isomorphic to a quadric
 in $\PP^3_k$.
\end{Proposition}

\prf
To show $(1)\Rightarrow(2)$, assume that ${\rm Deg}_X$ has a $k$-rational point.
Thus, there exists degenerate quadric  $Q$ with $X\subset Q\subset\PP^4_k$.
As explained in the proof of \cite[Theorem 8.6.8]{Dolgachev}, $Q$ is a cone over a smooth 
quadric surface, and the projection away from its vertex $\PP^4_k\dashrightarrow\PP^3_k$ induces
a morphism $X\to \PP^3_k$ that is finite of degree $2$ onto a smooth quadric surface $S$.
In particular, $S$ is a del~Pezzo surface of degree $8$ of product type.

To show $(2)\Rightarrow(1)$, let $\psi:X\to S$ be as in the statement.
Then, we have a short exact sequence (which even splits since ${\rm char}(k)\neq2$)
$$ 
  0\,\to\,\OO_S\,\to\,\psi_\ast\OO_X\,\to\,{\calL}^{-1}\,\to\,0,
$$
where $\calL$ is an invertible sheaf on $S$, which is of type $(1,1)$ 
on $S_{\overline{k}}\iso\PP^1_{\overline{k}}\times\PP^1_{\overline{k}}$.
In particular, $|\calL|$ defines an embedding $\imath:S\to\PP^3_k$ as a quadric,
and establishes the final assertion.
Now, $\imath\circ\psi$ arises from a $4$-dimensional subspace $V$ inside the
linear system $(\imath\circ\psi)^*\OO_{\PP^3_k}(1)\iso\omega_X^{-1}$.
Thus, $\imath\circ\psi$ is the composition of the anti-canonical embedding 
$X\to\PP^4_k$, followed by a projection $\PP^4_k\dashrightarrow\PP^3_k$.
As explained in the proof of \cite[Theorem 8.6.8]{Dolgachev}, such a projection
induces a degree $2$ morphism onto a quadric if and only if the point of projection
is the vertex of a singular quadric in $\PP^4_k$ containing $X$.
In particular, this vertex and the corresponding quadric are defined over $k$,
giving rise to a $k$-rational point of ${\rm Deg}_X$.
\qed\medskip

In order to refine Proposition \ref{prop: dP4 to dP8}, we will use
conic classes as introduced in Definition \ref{def: conic}.

\begin{Proposition}
  Let $X$ be a del~Pezzo surface of degree $4$ over a field $k$.
 Then, the following are equivalent:
  \begin{enumerate}
   \item There exists a conic class in $\Pic_{(X/k)({\rm fppf})}(k)$.
   \item There exists a finite morphism $\psi:X\to P'\times P''$ of degree $2$, where 
    $P'$ and $P''$ are a Brauer--Severi curves over $k$.
  \end{enumerate}
   Moreover, if $\psi$ is as in $(2)$, then $P'\iso P''$.
\end{Proposition}

\prf
Let $\calL\in\Pic_{(X/k)({\rm fppf})}(k)$ be a conic class.
By Theorem \ref{thm: main}, there exist morphisms 
$|{\calL}|:X\to P'$ and $|\omega_X^{-1}\otimes\calL^{-1}|:X\to P''$, where
$P'$ and $P''$ are Brauer--Severi curves over $k$.
Combining them, we obtain a finite morphism $X\to P'\times P''$ of degree $2$.
As in the proof of $(2)\Rightarrow(1)$ of Proposition \ref{prop: dP4 to dP8} 
we find that $P'\times P''$ embeds into $\PP^3$,
and thus, $0=[\PP^3_k]=[P']+[P'']\in\Br(k)$ by Proposition \ref{prop: product type classification}.
This implies $[P']=[P'']$ since these classes are $2$-torsion, and
thus, $P'\iso P''$ by Corollary \ref{cor: isomorphic BS}.
This establishes $(1)\Rightarrow(2)$.

Conversely, let $\psi:X\to P'\times P''$ be as in (2).
Then, $\psi^*(\OO_{P'}(1)\boxtimes\OO_{P''}(1))$ is a conic class, 
and (1) follows.
\qed\medskip

\subsection{Del~Pezzo surfaces of degree 3}
For these surfaces, we have the following analogue of 
Theorem \ref{thm: del Pezzo descent}.

\begin{Theorem}
  Let $X$ be a del~Pezzo surface of degree $3$ over a field $k$.
  Then, the following are equivalent:
  \begin{enumerate}
   \item There exists a birational morphism $f:X\to P$ to a Brauer--Severi surface,
     such that $f_{\overline{k}}$ is the blow-up of $6$ points in general position.
   \item There exists a class $F\in\Pic_{(X/k)({\rm fppf})}(k)$ such that 
     $$
        F_{\overline{k}} \,\iso\, \OO_{\overline{X}}(E_1+E_2+E_3+E_4+E_5+E_6),
     $$
     where the $E_i$ are disjoint $(-1)$-curves on $\overline{X}$.
 \end{enumerate}
\end{Theorem}

\prf
The proof is analogous to that of Theorem \ref{thm: dP5 descent},
and we leave the details to the reader.
\qed\medskip

Note that if the equivalent conditions of this theorem are fulfilled,
then $X$ is not minimal.
But the converse does not hold in general:
If $Y$ is a unirational, but not rational del~Pezzo surface of degree $4$ over
$k$, and $y\in Y$ is a $k$-rational point not lying on an exceptional curve,
then the blow-up $X\to Y$ in $y$ is a non-minimal
degree $3$ del~Pezzo surface over $k$ with $k$-rational points
that is not birationally equivalent to a Brauer--Severi surface over $k$.

By \cite[Theorem 28.1]{Manin}, a degree $3$ del~Pezzo surface $X$
is minimal if and only if $\rho(X)=1$, i.e.,
$\Pic_{(X/k)({\rm fppf})}(k)=\ZZ\cdot\omega_X$.
In this case, we have ${\rm Am}(X)=0$.
In particular, if such a surface is birationally equivalent to a Brauer--Severi surface $P$,
then $P\iso\PP^2_k$ by Proposition \ref{prop: amitsur is birational invariant}
and Theorem \ref{thm: amitsur birational}.
In particular, $X$ is rational and has a $k$-rational point in 
this case.

\subsection{Del~Pezzo surfaces of degree 2}
Arguing as in the proof of Theorem \ref{thm: birational to BS surface},
it follows that if there exists a dominant and rational map $P\dashrightarrow X$ from a 
Brauer--Severi surface $P$ onto a degree $2$ del~Pezzo surface over a perfect field $k$,
then $X$ has a $k$-rational point, and thus ${\rm Am}(X)=0$.
In particular, if $X$ is birationally equivalent to a Brauer--Severi surface, then it is rational, see
also Corollary \ref{cor: dP has a rational point}.

By work of Manin \cite[Theorem 29.4]{Manin}, a degree $2$ del~Pezzo surface over a field $k$
is unirational if it has a $k$-rational point not lying on an exceptional curve.
Together with non-trivial refinements of \cite{Salgado} and \cite{Festi}, such surfaces
over finite fields are always unirational.

By Lemma \ref{lem: del Pezzo Amitsur}, we have that ${\rm Am}(X)$ is at most
$2$-torsion for degree $2$ del~Pezzo surfaces.
For the possibilities of $H^1(G_k,\Pic_{X/k}(k^{\rm sep}))$, as well as further information
concerning arithmetic questions, we refer to \cite{Kresch Tschinkel}.

\subsection{Del~Pezzo surfaces of degree 1}

If $X$ is a del~Pezzo surface of degree $1$, then it has a $k$-rational point,
namely the unique base point of $|-K_X|$.
Thus, we have ${\rm Am}(X)=0$, and there are no morphisms or rational
maps to non-trivial Brauer--Severi varieties.

\end{document}